\documentclass[11pt,reqno]{amsart}

\usepackage{amssymb}
\usepackage{geompsfi}

\catcode`\@=11

\long\def\@savemarbox#1#2{\global\setbox#1\vtop{\hsize\marginparwidth 
  \@parboxrestore\tiny\raggedright #2}}
\newcommand\lref[1]{\ref{#1}%
\@ifundefined{r@DisplaY #1}{}{ (#1)}}

\newcommand\fakelabel[2]{\@bsphack\if@filesw {\let\thepage\relax
   \newcommand\protect{\noexpand\noexpand\noexpand}%
\xdef\@gtempa{\write\@auxout{\string
      \newlabel{#1}{{#2}{\thepage}}}}}\@gtempa
   \if@nobreak \ifvmode\nobreak\fi\fi\fi\@esphack}

\catcode`\@=12

\def\Empty{}
\newcommand\oplabel[1]{
  \def\OpArg{#1} \ifx \OpArg\Empty {} \else
        \label{#1}
  \fi}
\newtheorem{theoremSt}{Theorem}[section]

\newtheorem{exampleSt}[theoremSt]{Example}
\newtheorem{exerciseSt}[theoremSt]{Exercise}

\newcommand\MakeStEnv[1]{
  \newenvironment{#1}[1]{
  \begin{#1St} \oplabel{##1}%
  \global\def\CrntSt{\thetheoremSt}%
}{ 
  \end{#1St} }
  \newenvironment{#1+}[1]{
  \begin{#1St} \label{##1}%
  \label{DisplaY ##1}%
  \global\def\CrntSt{\thetheoremSt}%
  \def\Labl{##1}\ifx\Labl\Empty{} \else {\em (\Labl)\,}\fi%
}{ 
  \end{#1St} }
}
\MakeStEnv{theorem}
\MakeStEnv{corollary}
\MakeStEnv{proposition}
\MakeStEnv{lemma}
\MakeStEnv{definition}
\MakeStEnv{conjecture}
\MakeStEnv{problem}
\MakeStEnv{question}

\long\def\state#1#2{
\medskip\par\noindent
{\bf #1} 
{\it #2}
\par\medskip
}

\long\def\realfig#1#2{
\begin{figure}[htbp]
\centerline{\psfig{file=#1}}
\caption[#1]{#2}
\oplabel{#1}
\end{figure}}

\newlength{\saveu}

\newenvironment{pf*}[1]{%
 \begin{proof}[#1]%
}{ 
 \end{proof}
}

\newcommand{\finishproof}[1]{ 
  \def\FPArg{#1}
  \ifx\FPArg\Empty
        \newcommand\FPArg{\CrntSt}  \fi
  \smallbreak\noindent\makebox[\textwidth]{\hfill\fbox{\FPArg}}
  \medbreak\noindent
}

\newcommand\AAA{{\mathcal A}}

\newcommand\CC{{\mathcal C}}

\newcommand\FF{{\mathcal F}}

\newcommand\LL{{\mathcal L}}
\newcommand\MM{{\mathcal M}}

\newcommand\PP{{\mathcal P}}

\newcommand\UU{{\mathcal U}}

\newcommand\PMF{{\PP\kern-2pt\MM\FF}}
\newcommand\ML{{\MM\LL}}
\newcommand\PML{{\PP\kern-2pt\MM\LL}}

\newcommand\half{{\textstyle{\frac12}}}

\newcommand\ep{\epsilon}

\newcommand\hhat{\widehat}

\newcommand\union{\cup}
\newcommand\intersect{\cap}
\newcommand\bbR{{\mathord{\text{I\kern-2pt R}}}}        
\newcommand\bbH{{\mathord{\text{I\kern-2pt H}}}}        

\newcommand\C{{\mathbb C}}

\newcommand\Z{{\mathbb Z}}
\newcommand\R{{\mathbb R}}
\newcommand\Q{{\mathbb Q}}

\newcommand\Hyp{{\mathbb H}}

\newcommand\PSL[1]{\text{PSL}_{#1}}

\newcommand\bigrightarrow[1]{\hbox to #1{\rightarrowfill}}
\newcommand\bigleftarrow[1]{\hbox to #1{\leftarrowfill}}

\newcommand\homeo{\cong}
\newcommand\boundary{\partial}
\newcommand\semidir{\mathrel{\hbox{\vrule depth-.03ex height1.1ex\kern-0.15em$\times$}}}

\newcommand{\diam}{\operatorname{diam}}

\renewcommand{\Re}{\operatorname{Re}}
\renewcommand{\Im}{\operatorname{Im}}

\numberwithin{equation}{section}

\catcode`\@=11

\def\subsection{\@startsection{subsection}{2}%
  \z@{.5\linespacing\@plus.7\linespacing}{.5em}%
  {\normalfont\bfseries\centering}}

\def\section{\@startsection{section}{1}%
  \z@{.7\linespacing\@plus\linespacing}{.5\linespacing}%
  {\normalfont\large\bfseries\centering}}

\def\subsubsection{\@startsection{subsubsection}{3}%
  \z@{.5\linespacing\@plus.7\linespacing}{-.5em}%
  {\normalfont\bfseries}}

\catcode`\@=12

\newcommand{\pleat}{\operatorname{\mathbf{pleat}}}
\newcommand{\short}{\operatorname{\mathbf{short}}}

\newcommand{\dist}{\operatorname{dist}}

\newcommand{\fsubd}{\mathrel{{\scriptstyle\searrow}\kern-1ex^d\kern0.5ex}}
\newcommand{\bsubd}{\mathrel{{\scriptstyle\swarrow}\kern-1.6ex^d\kern0.8ex}}
\newcommand{\fsubeq}{\mathrel{\genfrac{}{}{0pt}{}{\searrow}{\raise.1ex\hbox{=}}}}
\newcommand{\bsubeq}{\mathrel{\genfrac{}{}{0pt}{}{\swarrow}{\raise.1ex\hbox{=}}}}

\newcommand{\bbar}{\overline}
\newcommand{\UML}{\operatorname{\UU\MM\LL}}
\newcommand{\EL}{\mathcal{EL}}

\newcommand{\tsh}[1]{\left\{\kern-.9ex\left\{#1\right\}\kern-.9ex\right\}}

\newcommand\Momega{\omega_M}

\newcommand\modl{M_\nu}
\newcommand\MT{{\mathbb T}}

\newcommand\nbhd{\operatorname{Nbhd}}

\begin{document}

\title[Aspects of Hyperbolic 3-manifolds]{Combinatorial and
Geometrical Aspects of Hyperbolic 3-Manifolds} 

\author{Yair N. Minsky}
\address{SUNY Stony Brook}
\date{\today}


\maketitle


 
\section{Introduction}

This is the edited and revised form of handwritten
notes that were distributed with the lectures that
I gave at the workshop on Kleinian Groups and Hyperbolic 3-Manifolds
in Warwick on September 11-15 of 2001.\footnote{The 
terrible events in New York that coincided with the beginning of this
conference overshadow its subject matter in significance,
and yet those same events demand of us to continue with our ordinary work.}
The goal of the lectures was
to expose some recent work \cite{minsky:ELCI}
on the structure of ends of hyperbolic
3-manifolds, which 
is part of a program to solve Thurston's Ending
Lamination Conjecture (the conclusion of the program, which is joint
work with J. Brock and R. Canary, will appear in
\cite{brock-canary-minsky:ELCII}). In the interests of simplicity and
the ability 
to get to the heart of the matter, the notes are quite informal in
their treatment of background material, and the main results are often
stated in special cases, with detailed examples taking the place of
proofs. Thus it is hoped that the reader will be able to extract the
main ideas with a minimal investment of effort, and in the event he or
she is still interested, can obtain the details in \cite{minsky:ELCI},
which will appear later on.

I would like to thank the organizers of the conference for inviting me
and giving me the opportunity to talk for what must have seemed like a
very long time. 

\subsection*{Object of Study} 

If the interior $N$ of a compact 3-manifold $\bbar N$ admits a
complete infinite-volume hyperbolic structure, then there is a
multidimensional {\em deformation space} of such structures. The study
of this space goes back to Poincar\'e and Klein, but the modern theory
began with Ahlfors-Bers in the 1960's and
received the perspective that we will focus on from Thurston
and others in the late 70's. The deformation theory depends
deeply on an understanding of the geometry of the {\em ends} of $N$
(in the sense of Freudenthal
\cite{freudenthal:ends}), which one can think of as small
neighborhoods of the boundary components of $\bbar N$.

The interior of the deformation space, as studied by 
Ahlfors, Bers
\cite{ahlfors-bers,bers:simultaneous,bers:spaces},
Kra \cite{kra:spaces}, Marden
\cite{marden-maskit,marden:geometry},
Maskit \cite{maskit:classification},  and Sullivan 
\cite{sullivan:QCDII},
can be parametrized using the Teichm\"uller space of 
$\boundary \bbar N$ -- that is, by choosing a ``conformal structure at
infinity'' for each (non-toroidal) boundary component of $\bbar N$. 
(See also \cite{keen-series:pleating} and \cite{bonahon-otal:bending}
for other approaches to the study of the interior).
The boundary contains manifolds with parabolic cusps
\cite{maskit:boundaries,mcmullen:cusps}, 
and more generally,
with {\em geometrically infinite ends}
\cite{bers:boundaries,greenberg:geominf,wpt:notes}. 
The Teichm\"uller parameter is replaced by Thurston's {\em
ending laminations} for such ends. Thurston conjectured
\cite{wpt:bull} that these invariants are sufficient to determine the
geometry of $N$ uniquely -- this is known as the Ending Lamination
Conjecture (see also \cite{abikoff:survey} for a survey). 

In these notes we will consider the special case of {\em Kleinian
surface groups}, for which $\pi_1(N)$ is isomorphic to $\pi_1(S)$ for
a surface $S$. This case suffices for describing the ends of general
$N$, provided $\boundary \bbar N$ is incompressible.
(In the compressible case the deeper question of Marden's {\em
tameness conjecture} comes in, and this is beyond the scope of our
discussion. See  Marden \cite{marden:geometry} and Canary
\cite{canary:ends}.)

We will show how the the ending laminations, together
with the combinatorial structure of the set of simple closed curves on
a surface,  allows us to
build a {\em Lipschitz model} for the geometric structure of $N$,
which in particular describes the thick-thin decomposition of
$N$. These results, which are proven in detail in \cite{minsky:ELCI}, 
will later be followed by {\em bilipschitz} estimates
in Brock-Canary-Minsky \cite{brock-canary-minsky:ELCII}, 
and these will suffice to prove Thurston's conjecture in the case of
incompressible boundary.

\subsection*{Kleinian surface groups}

From now on, let
$S$ be an oriented compact surface
with $\chi(S)<0$, and let
$$\rho:\pi_1(S) \to \PSL 2(\C)$$ 
be a discrete, faithful representation.
If $\boundary S \ne \emptyset$ we require $\rho(\gamma)$ to be
parabolic for $\gamma$ representing any boundary component. 
This is known as a (marked) Kleinian surface group.
We name the quotient 3-manifold
$$N=N_\rho = \Hyp^3/\rho(\pi_1(S)).$$

\subsubsection*{Periodic manifolds}
Before discussing the general situation let us consider a well-known
and especially tractable example. 

Let $\varphi:S\to S$ be a pseudo-Anosov homeomorphism (this means that
$\varphi$ leaves no finite set of non-boundary curves invariant up to
isotopy). The mapping torus of $\varphi$ is
$$
M_\varphi = S\times \R / \langle (x,t)\mapsto (\varphi(x),t+1)\rangle,
$$
a surface bundle over $S^1$ with fibre $S$ and  monodromy $\varphi$. 
Thurston \cite{wpt:II} showed, as part of his hyperbolization theorem,
that $int(M_\varphi)$
admits a hyperbolic structure which we'll call $N_\varphi$ (see also Otal
\cite{otal:surfacebundles} and McMullen \cite{mcmullen:renormbook}).
Let  $N\homeo int(S)\times\R$ be the infinite cyclic cover of $N_\varphi$,
``unwrapping'' 
the circle direction (Figure \ref{bundlecover}).
After identifying $S$ with some lift of the fibre, 
we obtain an isomorphism $\rho:\pi_1(S)
 \to \pi_1(N)\subset \PSL 2(\C)$, which is a Kleinian surface group. 

\realfig{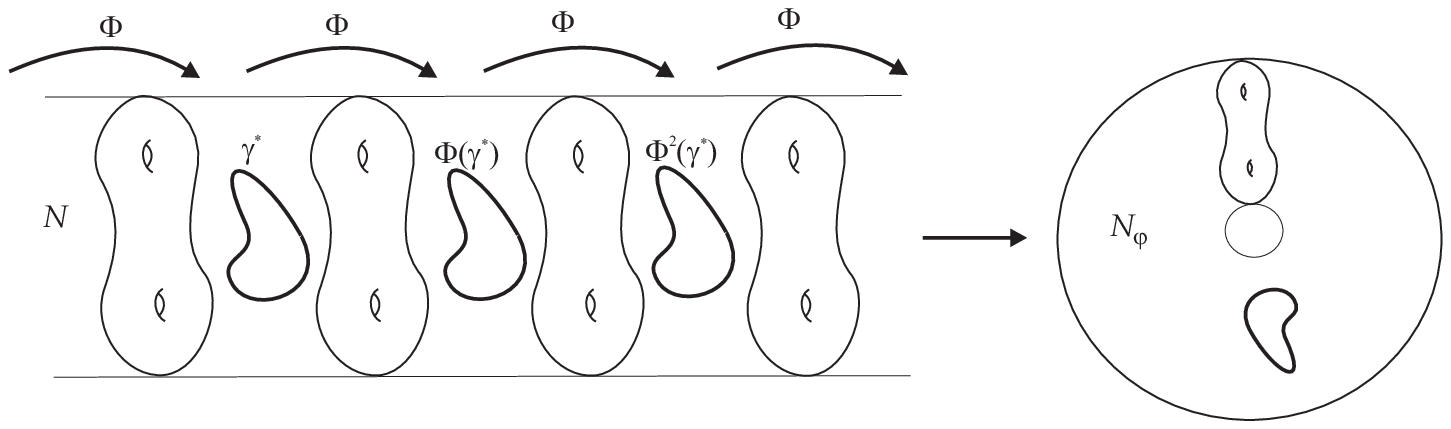}{$N$ covers the surface
bundle $N_\varphi$.}

The deck translation $\Phi:N\to N$ of the covering induces
$\Phi_* = \varphi_* : \pi_1(S) \to \pi_1(S)$. We next 
consider the action of $\varphi$
on the space of {\em projective measured laminations} $\PML(S)$
(see \cite{travaux,bonahon:laminations}, and Lecture 3). 
For every simple closed curve $\gamma$ in $S$, the sequences
$[\varphi^n(\gamma)]$ and $[\varphi^{-n}(\gamma)]$ converge to two
distinct points $\nu_+$ and $\nu_-$ in $\PML(S)$. After isotopy,
$\varphi$ can be represented on $S$ by a map that preserves the leaves
of both $\nu_+$ and $\nu_-$, stretching the former and contracting the
latter. 

We can see $\nu_\pm$ directly in the asymptotic geometry  of $N$:
For a curve $\gamma$ in $S$, let $\gamma^*$ be its geodesic
representative in $N$.
Now consider $\Phi^n(\gamma^*)$ -- these are all geodesics of the same
length, marching off to infinity in both directions as $n\to
\pm\infty$, 
and note that $\Phi^n(\gamma^*) = \varphi^n(\gamma)^*$.
So, we have a sequence of simple curves in $S$, converging to $\nu_+$ as
$n\to\infty$, whose geodesic representatives ``exit the $+$ end'' of
$N$ (similarly as $n\to\infty$ they converge to $\nu_-$ and the
geodesics exit the other end).

The laminations $\nu_\pm$ are the ending laminations of $\rho$ in this
case. To understand the general case we will have to develop a bit of
terminology, and recall the work of Thurston and Bonahon. 

\subsubsection*{Ends}
Let $N_0$ denote $N$ minus its cusps (each cusp is an open solid
torus, whose boundary in $N$ is a properly embedded open annulus). 
McCullough's relative version \cite{mccullough:relative} of Scott's 
core theorem \cite{scott:core} gives us a
compact submanifold $K$ in $N_0$, homeomorphic to $S\times[0,1]$,
which meets each cusp boundary in an annulus (including the annuli
$\boundary S \times[0,1]$). The components of $N_0 \setminus K$ are
in one to one correspondence with the topological ends of $N_0$, and
are called neighborhoods of the ends (see Bonahon \cite{bonahon}).

$N$ also has a {\em convex core} $C_N$, which is the smallest closed
convex submanifold whose inclusion is a homotopy equivalence. Each end
neighborhood either meets $C_N$ in a bounded set, in which case the
end is called {\em geometrically finite}, or is contained in $C_N$, in
which case the end is {\em geometrically infinite}.

From now on, let us assume that $N$ has {\em no extra cusps},
which means that the cusps correspond only to the components of
$\boundary S$. In particular $N_0$ has exactly two ends, which we
label $+$ and $-$ according to an appropriate convention.

\subsubsection*{Simply degenerate ends}
In \cite{wpt:notes}, Thurston made the following definition, which can
be motivated by the surface bundle example: 

\begin{definition}{def degenerate}
An end of $N$ is {\bf simply degenerate} if there exists a sequence of
simple closed curves $\alpha_i$ in $S$ such that $\alpha_i^*$ exit the
end.
\end{definition}

Here ``exiting the end'' means that the geodesics are eventually
contained in an arbitrarily small neighborhood of the end, and in
particular outside any compact set. 
 Note that a geometrically finite
end cannot be simply degenerate, since all closed geodesics are
contained in the convex hull. 

Thurston then established this theorem (stated in the case without
extra cusps):

\begin{theorem}{EL exist}{\rm \cite{wpt:notes}}
If an end $e$ of $N$ is simply degenerate then there exists a unique
lamination $\nu_e$ in $S$ such that for any sequence of simple closed
curves $\alpha_i$ in $S$, 
$$
\alpha_i \to \nu_e \iff \text{$\alpha_i^*$ exit the end $e$.}
$$
A sequence $\alpha_i\to \nu_e$ can be chosen so that the lengths
$\ell_N(\alpha_i^*) \le L_0$, where $L_0$ depends only on $S$.

Furthermore, $\nu_e$ {\em fills} $S$ -- its complement consists of
ideal polygons and once-punctured ideal polygons.
\end{theorem}

(We are being cagey here about just what kind of lamination $\nu_e$ is,
and what convergence $\alpha_i\to\nu_e$ means. See Lecture 3
for more details.)

Thurston also proved that simply degenerate ends are
{\em tame}, meaning that they have neighborhoods homeomorphic to
$S\times(0,\infty)$,
and that manifolds obtained as limits of quasifuchsian manifolds
have ends that are geometrically finite or simply degenerate. Bonahon
completed the picture with his ``tameness theorem'', 

\begin{theorem}{Bonahon tameness}{\rm \cite{bonahon}}
The ends of $N$ are either geometrically finite or simply degenerate.
\end{theorem}

In particular $N_0$ is homeomorphic to $S\times\R$, and ending
laminations are well-defined for each geometrically infinite end.

Geometrically finite ends are the ones treated by  Ahlfors,
Bers and their coworkers, and their analysis requires a discussion of
quasiconformal mappings and Teichm\"uller theory
(see \cite{bers:simultaneous,bers:spaces,sullivan:ICM} for more). In
order to simplify our exposition
we will limit ourselves, for the remainder of these
notes, to Kleinian surface groups $\rho$ with no extra cusps, 
and with no geometrically finite ends. In particular the convex hull
of $N_\rho$ is all of $N_\rho$, and there are two ending laminations,
$\nu_+$ and $\nu_-$.  This is called the
{\em doubly degenerate} case.

\subsection*{Models and bounds}
Our goal now is to recover geometric information about $N_\rho$ from
the asymptotic data encoded in $\nu_\pm$.  The following natural
questions arise, for example: 

\begin{itemize}
\item
Thurston's Theorem \ref{EL exist} guarantees the existence of a sequence 
$\alpha_i \to \nu_+$ whose geodesic representatives have bounded lengths
$\ell_N(\alpha_i^*)$. How can we determine, from $\nu_+$,
which sequences have this property?
\item
The case of the cyclic cover of a
surface bundle is not typical: because it covers a compact manifold
(except for cusps), it
has ``bounded geometry''. That is, 
$$\inf_\beta \ell_N(\beta) > 0$$ 
where $\beta$ varies over closed geodesics. The bounded geometry case
is considerably easier to understand. In particular the Ending
Lamination Conjecture in this category (without cusps) was proven in
\cite{minsky:slowmaps,minsky:endinglam}.

\noindent
Can we tell from $\nu_\pm$ alone whether $N$ has bounded geometry?

\item
If $N$ doesn't have bounded geometry, there are arbitrarily short
closed geodesics in $N$, each one encased in a {\em Margulis tube},
which is a standard collar neighborhood. Such examples were shown to
exist by Thurston \cite{wpt:II} and Bonahon-Otal \cite{bonahon-otal},
and to be generic in an appropriate sense by McMullen \cite{mcmullen:cusps}.

\noindent
In the unbounded geometry case, 
can we tell {\em which} curves
in $N$ are short? How are they arranged in $N$?
\end{itemize}

We will describe the construction of a ``model manifold'' $\modl$ for
$N$, which can be used to answer these questions. $\modl$ is
constructed combinatorially from $\nu_\pm$, and contains for example
solid tori that correspond to the Margulis tubes of short curves in
$N$. $\modl$ comes equipped with a map
$$
f:\modl \to N
$$
which takes the solid tori to the Margulis tubes, is proper, 
Lipschitz in the complement of the solid tori, and
preserves the end structure. This will be the content of
the Lipschitz Model Theorem, which will be stated precisely in Lecture 6.

Note that if $f$ is {\em bilipschitz} then the Ending Lamination
Conjecture follows: If $N_1,N_2$ have the same invariants $\nu_\pm$
then the same model $\modl$ would admit bilipschitz maps $f_1:\modl
\to N_1$ and $f_2:\modl \to N_2$, and $f_2\circ f_1^{-1} : N_1 \to
N_2$ would be a bilipschitz homeomorphism. By Sullivan's Rigidity
Theorem \cite{sullivan:rigidity}, $N_1$ and $N_2$ would be isometric.

\subsection*{Plan}

Here is a rough outline of the remaining lectures:
\begin{description}
\item[\S 2] {\bf Hierarchies and model manifolds:} 
We will show how to build $\modl$ starting with a geodesic in the
{\em complex of curves} $\CC(S)$. The main tool is the {\em hierarchy
of geodesics} developed in Masur-Minsky \cite{masur-minsky:complex2}.
Much of the discussion will take place in the special case of the
5-holed sphere $S_{0,5}$, where the definitions and arguments are
considerably  simplified.

\item[\S 3] {\bf Ending laminations to model:}
Using a theorem of Klarreich we will relate ending laminations to
{\em points at infinity} for $\CC(S)$, and this will allow us to
associate to a pair of ending laminations a geodesic in $\CC(S)$, and
its associated hierarchy and model manifold. 

\item[\S 4] {\bf Quasiconvexity:}
We then begin to explore the linkage
between geometry of the
3-manifold $N_\rho$ and the curve complex data. We will show that the
subset of $\CC(S)$  consisting of curves with bounded length in $N$
is {\em quasiconvex}.
The main tool here is an argument using pleated surfaces and
Thurston's Uniform Injectivity Theorem. 

\item[\S 5] {\bf Projection Bounds:}
In this lecture we will discuss the
Projection Bound Theorem, a
strengthening of the Quasiconvexity Theorem
that shows that curves that
appear in the hierarchy are combinatorially close
to the bounded-length curves in $N$. We will also prove the Tube
Penetration Theorem, which controls how deeply certain pleated
surfaces can enter into Margulis tubes. 

\item[\S 6]{\bf A priori bounds and the model map:}
Applying the Projection Bound Theorem and the Tube Penetration
Theorem, we will establish a uniform bound on the lengths of all
curves that appear in the hierarchy.

We will then state the Lipschitz Model Theorem, whose proof uses the
a priori bound and a few additional geometric arguments. As
consequences we will obtain some final statements on the structure of
the set of short curves in $N$.

\end{description}

\section{Curve Complex and Model Manifold}

In this lecture we will introduce the complex of curves $\CC(S)$ and
demonstrate how a geodesic in $\CC(S)$ leads us to construct a ``model
manifold''. For simplicity we will mostly work with $S=S_{0,5}$, the
sphere with 5 holes.  (In general let $S_{g,n}$ be the surface with
genus $g$ and $n$ boundary components).

\subsection*{The complex of curves}
$\CC(S)$ will be a simplicial complex whose vertices are homotopy
classes of simple, essential, unoriented closed curves (``Essential''
means homotopically nontrivial, and not homotopic to the boundary).
Barring the exceptions below, we define 
the $k$-simplices to be unordered $k+1$-tuples $[v_0\ldots v_k]$
such that $\{v_i\}$ can be realized as pairwise disjoint curves. 
This definition was given by Harvey \cite{harvey:boundary}. 

\subsubsection*{Exceptions:} If $S=S_{0,4}$, $S_{1,0}$ or $S_{1,1}$
then this definition gives no edges. Instead we allow edges $[vw]$ whenver
$v$ and $w$ can be realized with 
$$
\#v\intersect w = \left\{\begin{array}{ll}
			1 & S_{1,0}, S_{1,1} \\
			2 & S_{0,4}.
		\end{array}\right.
$$
(see Figure \ref{fareycases}).
In this case $\CC(S)$ is the {\em Farey graph} in the plane: a vertex
is indexed by the slope $p/q$ of its lift to the planar $\Z^2$ cover of $S$,
so the vertex set is $\hhat \Q = \Q\union\infty$. 
Two vertices $p/q$, $r/s$ are joined by an edge if $|ps-qr|=1$ 
(see e.g. Series \cite{series:markoff} or \cite{minsky:torus}).

\realfig{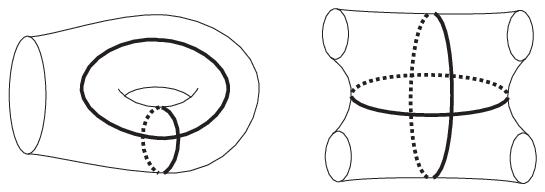}{Adjacent vertices in $\CC(S_{1,1})$ and $\CC(S_{0,4})$}

For $S_{0,0}, S_{0,1}, S_{0,2}, S_{0,3}$: $\CC(S)$ is empty.
(For the annulus $S_{0,2}$ there is another useful construction which
we will return to later.)

Let $\CC_k(S)$ denote the $k$-skeleton of $\CC(S)$. We will
concentrate on $\CC_0$ and $\CC_1$. 

We endow $\CC(S)$ with the metric  that makes every simplex regular
Euclidean of sidelength 1. Thus $\CC_1(S)$ is a graph with unit-length
edges. Consider a geodesic in $\CC_1(S)$ -- it is a sequence of
vertices $\{v_i\}$ connected by edges (Figure \ref{geodcases}),
and in particular: $v_i,v_{i+1}$ are disjoint (in the non-exceptional
cases), $v_i$ and $v_{i+2}$ intersect but are disjoint from $v_{i+1}$,
and $v_i$ and $v_{i+3}$ {\em fill} the surface: their union intersects
every essential curve. It is harder to characterize topologically the
relation between $v_i$ and $v_j$ for $j>i+3$.

\realfig{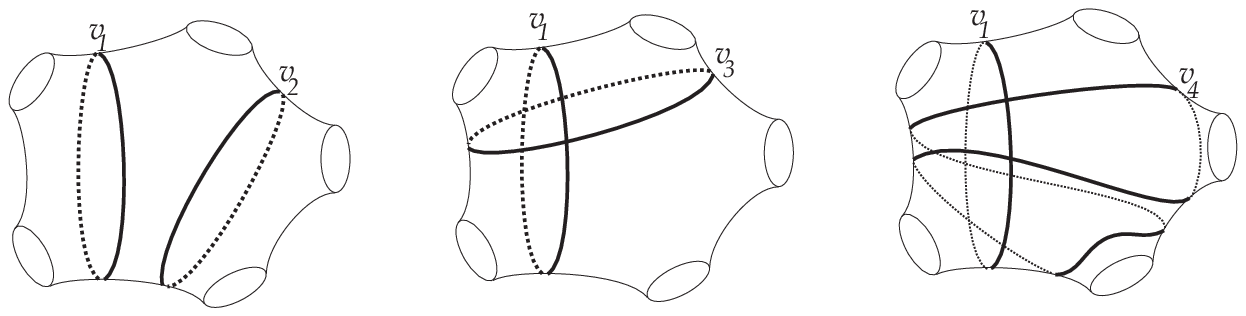}{$\{v_1,v_2,v_3,v_4\}$ are the vertices of a
geodesic in $\CC(S_{0,5})$.} 

\subsection*{Model Construction}

Let $S=S_{0,5}$ -- this case is considerably simpler than the general
case, while preserving many of the main features. 

Starting with a {\em bi-infinite geodesic} $g$ in $\CC_1(S)$ (more
about the existence of such geodesics later), we will construct a
manifold $M_g\homeo S\times\R$, equipped with a piecewise-Riemannian
metric. $M_g$ is made of ``standard blocks'', all isometric, and
``tubes'', or solid tori of the form
$\text{(annulus)}\times\text{(interval)}$.

\subsubsection*{Hierarchy}
We begin by ``thickening'' $g$ in the following sense: 
Any vertex $v\in\CC_0(S)$ divides $S$ into two components, one $S_{0,3}$ and
one $S_{0,4}$. Let $W_v$ denote the second of these. 
If $v_i$ is a vertex of $g$ then 
$$
v_{i-1},v_{i+1} \in \CC_0(W_{v_i}).
$$
The complex $\CC(W_{v_i})$ is just the Farey graph, and we may join
$v_{i-1}$ to $v_{i+1}$ by a geodesic in that graph. Name this geodesic
$h_i$, and represent it schematically as in Figure \ref{wheels}.

\realfig{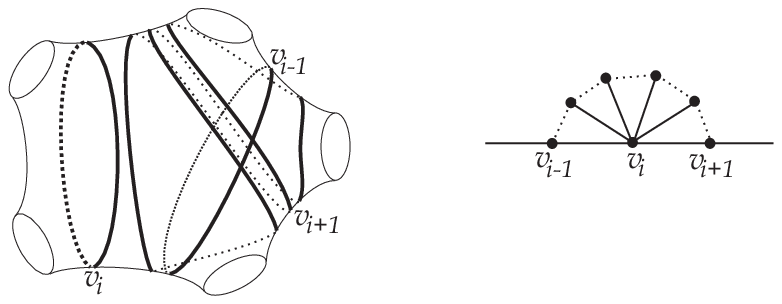}{The local configuration at
a vertex $v_i$ of $g$ yields a ``wheel'' in the link of $v_i$. Note,
edges of $h_i$ are not edges of $\CC(S)$; call them ``rim'' edges. The
other edges are called ``spokes''.}

We repeat this at every vertex. The resulting system is 
called a {\em hierarchy of geodesics}. 
(In general surfaces, considerable complications arise. Geodesics must
satisfy a technical condition called ``tightness'', and the hierarchy
has more levels. This is joint work with Masur
\cite{masur-minsky:complex1,masur-minsky:complex2}.)

Note that the construction is not uniquely dependent on $g$ -- there
are arbitrary choices for each $h_i$. However what we have to say will
work regardless of how the choices are made.

\subsubsection*{Blocks}
To each rim edge $e$ we associate a ``block'' $B(e)$, and then glue
these together to form the model manifold. 
$e$ is an edge of $\CC(W_v)$ for some vertex $v$ -- denote $W_e \equiv W_v$
for convenience. Let $e^-,e^+$ be its vertices, ordered from left to right.
Let $C_+$ and $C_-$ be open collar neighborhoods of $e^+$ and
$e^-$, respectively. We define
$$
B(e) = W_e \times[-1,1] - \left( C_+\times(1/2,1] \union C_- \times
[-1,-1/2) \right).
$$
Thus we have removed solid-torus ``trenches'' from the top and bottom
of the product $W_e\times[-1,1]$. Figure \ref{block} depicts this as
a gluing construction.

\realfig{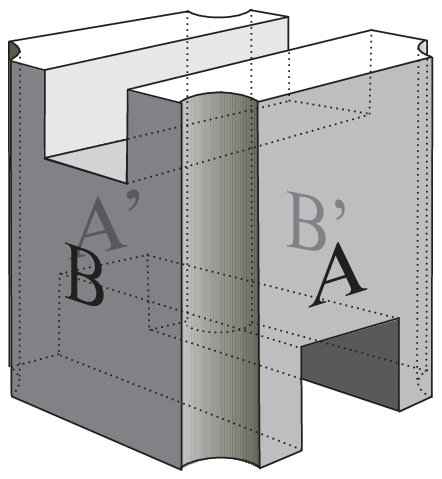}{Construct a block $B(e)$
by doubling this object along $A,A', B$ and $B'$. 
The curved vertical faces become
$\boundary W_e\times [-1,1]$.}

The boundary $\boundary B(e)$ divides into four {\em 3-holed spheres},
$$
\boundary_\pm B(e) \equiv (W_e-C_{\pm})\times{\pm 1}
$$
and some {\em annuli}. Schematically, we depict this structure in
Figure \ref{blockschem}.

\realfig{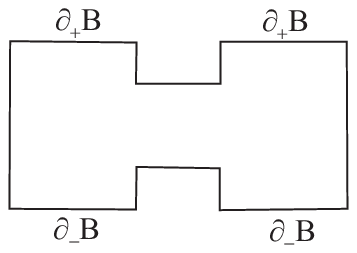}{Schematic diagram of the
  different pieces of the boundary of a block.}

\subsubsection*{Gluing}
Take the disjoint union of all the blocks arising from the hierarchy
over $g$, and glue them along 3-holed sphere, where possible. That is,
if $Y\times\{1\}$ appears in $\boundary_+B(e_1)$ and $Y\times\{-1\}$
appears in $\boundary_-B(e_2)$, identify them using the identity map in
$Y$. 

(A technicality we are eliding is that subsurfaces are determined only
up to isotopy; one can select one representative for each isotopy
class in a fairly nice and consistent way.)

There are three types of gluings that can occur:
\begin{enumerate}

\item Both edges occur in the link of the same vertex $v$; 
$W_{e_1} = W_{e_2} = W_v$, and $e_1^+ = e_2^-$ (Figure
\ref{gluing1}). $B(e_1)$ and $B(e_2)$ are glued along
$W_v\setminus C_{e_1^+}$, which is composed of three-holed sphere
$Y_1$ and $Y_2$. 

\realfig{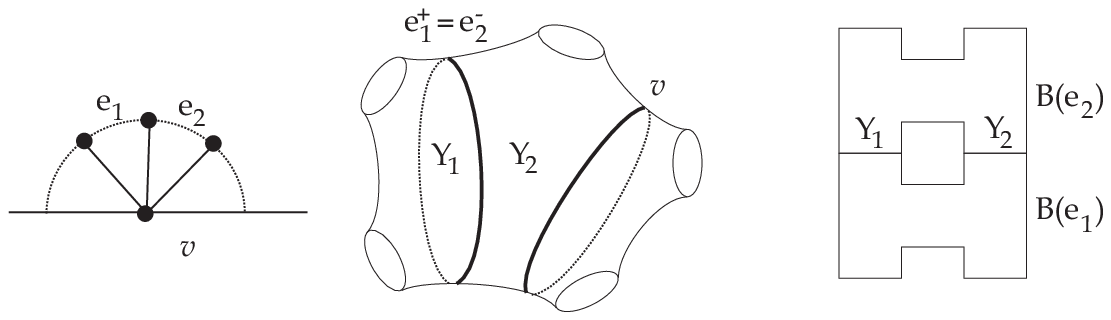}{}

\item $e_1\subset \CC(W_u)$ and $e_2\subset \CC(W_v)$, where $u$ and
$v$ are two succesive vertices (Figure \ref{gluing2}).
Now $e_1^+ = v$ and $e_2^-=u$, and
the gluing is along $Y_2 = W_u \intersect W_v$, which
separates $S_{0,5}$.

\realfig{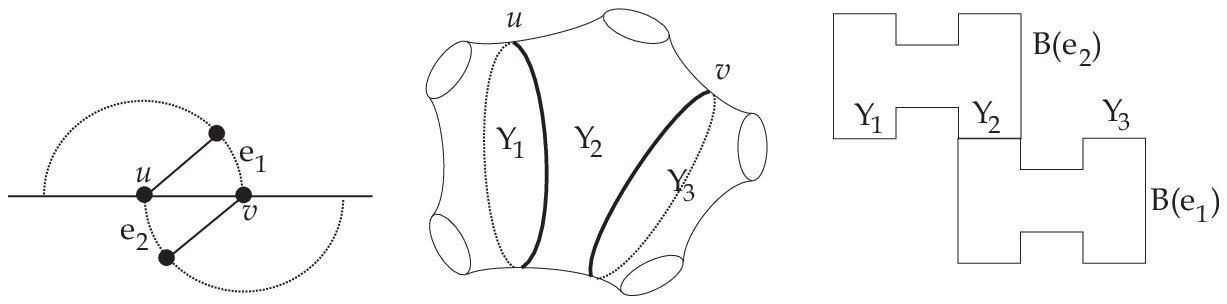}{}

\item $e_1\subset \CC(W_u)$ and $e_2\subset \CC(W_w)$, where $u,v,w$
are three successive vertices  (Figure \ref{gluing3}). 
In this case $e_1^+ = e_2^- = v$, and
the gluing is
along $Y_1$ which is isotopic to $W_u\intersect W_w$ and does not
separate $S_{0,5}$. Note that the intersection pattern of
$u$ and $w$ is typically more complicated than pictured, as
$d_{W_v}(u,w) >> 1$.

\realfig{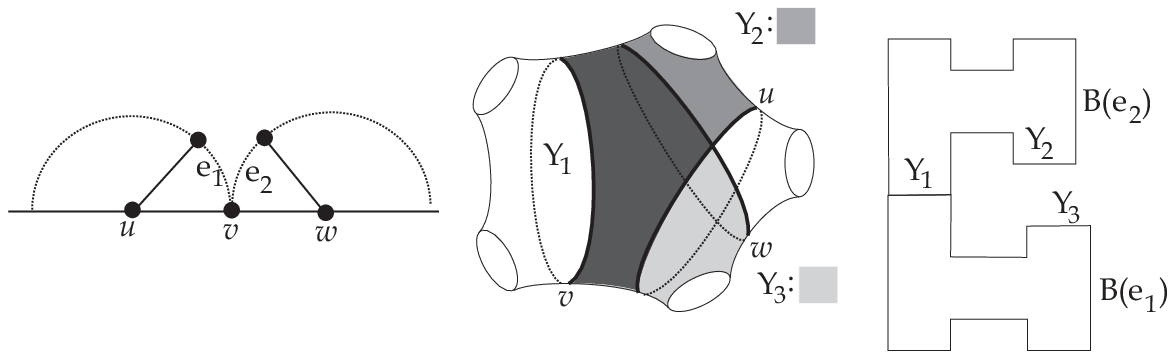}{}

\end{enumerate}

When we fit all the blocks together, the result can be embedded in
$S\times\R$, in such a way that any level surface $Y\times\{t\}$ in a
block is mapped to $Y\times\{s\}$ in $S\times\R$ by a map that is the
identity on the first factor. Call such a map ``straight''. Note that
the blocks can be stretched {\em vertically} in different ways.

\realfig{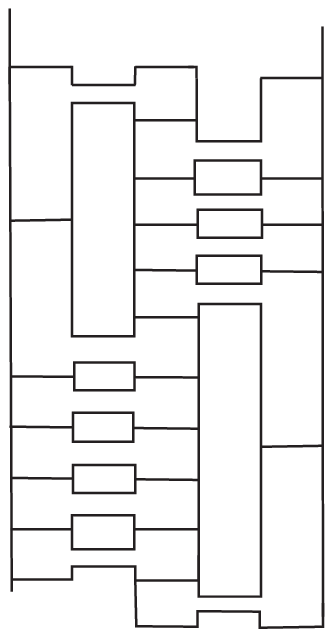}{A schematic of the embedding 
of the blocks in $S\times\R$}

In the gaps between blocks we find solid tori of the form
$$
C\times(s,t)
$$
where $C$ is one of our collar neighborhoods of a vertex in the
hierarchy. Call these the {\em  tubes} of the model. 

We should of course verify these claims about the gluing operations. 
Two things to check are: 
\begin{enumerate}
\item All the vertices in the hierarchy are distinct (and hence all
the tubes are homotopically distinct in $S\times\R$).
\item The gluings we have shown are the only ones.
\end{enumerate}

\realfig{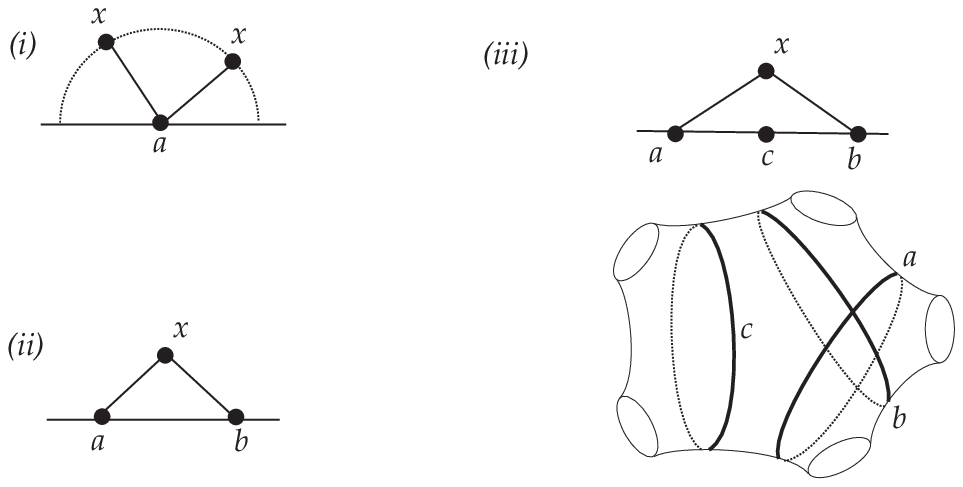}{The three ways that a vertex can
appear twice}

To verify (1), suppose that a vertex $x$ appears in two places in the
hierarchy. That is, $x$ is in the wheels
(links) of vertices $a$ and $b$ in $g$. The triangle inequality in
$\CC(S)$ implies that $d(a,b)\le 2$, and since $g$ is a geodesic this
leaves three possibilities, as in Figure \ref{vtwice}.

In case {\em (i)}, $a=b$. This is not possible since the ``rim'' path is a
geodesic in $\CC(W_a)$.

In case {\em (ii)}, $a,b$ and $x$ make a triangle in $\CC(S)$, but
$\CC(S)$ has no triangles for $S=S_{0,5}$. 

In case {\em (iii)}, $a$ and $b$ ``fill'' the 4-holed sphere $W_c$
bounded by $c$, so that if $x$ is
represented by a curve disjoint from both, it is equal to $c$ or lies
on the complement of $W_c$. That complement is a 3-holed sphere so the
only possibility is that $x=c$. In other words $x$ really only appears
once, as a vertex of $g$. 

\medskip

To prove (2), we must consider how a gluing surface $Y$ (a 3-holed
sphere) can occur. There are several possibilities for the curves of
$\boundary Y$  (Figure \ref{Yposs}).

\realfig{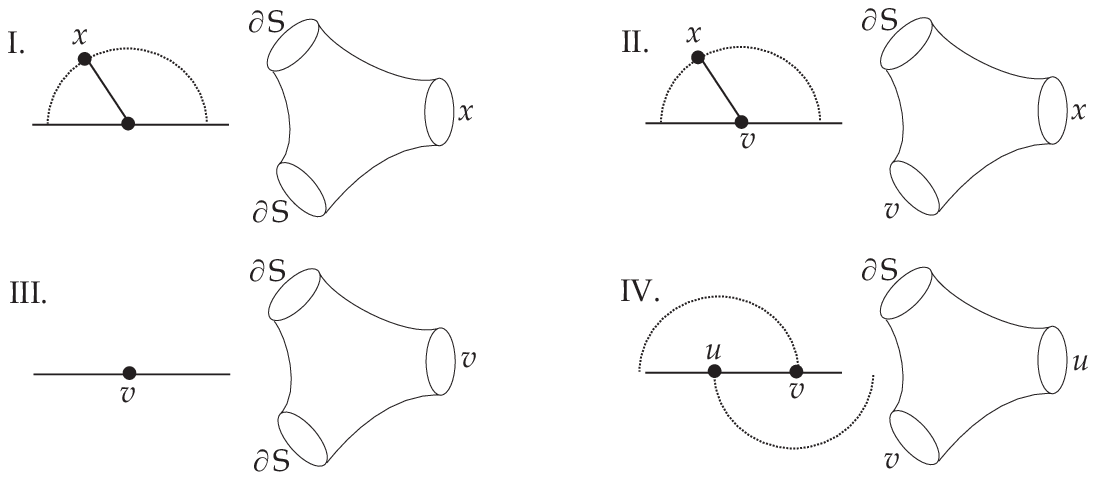}{}

\begin{enumerate}
\item[I.]  $\boundary Y$ consists of $x$ and two curves of $\boundary
S$, where $x$ is an interior ``rim'' vertex. 
\item[II.] $\boundary Y$ consists of an interior rim vertex $x$, a $g$ vertex
$v$, and a curve of $\boundary S$.
\item[III.] $\boundary Y$ consists of $v$ and two curves of $\boundary
S$, where $v$ is a vertex on $g$. 
\item[IV.] $\boundary Y$ consists of two adjacent vertices $u,v$ of
$g$ and a curve of $\boundary S$. 
\end{enumerate}
Types I and II occur in pairs as the top and bottom surfaces of two
blocks associated to adjacent rim edges meeting at the same $x$. 
Types III and IV occur on blocks associated to first and last edges in
rim geodesics, and each one occurs in exactly two ways. It is
therefore not hard to check all the possibilities and see that the
gluings we described indeed produce a manifold. 

The embedding of the manifold into $S\times\R$ can be done
inductively, by ``sweeping'' across the hierarchy from left to right.

\subsection*{Geometry of the model}

Fix one {\em standard block}: Take a copy of $W$ (of type $S_{0,4}$) with
two curves $\alpha,\beta$ that are neighbors in $\CC(W)$, collars
$C_\alpha,C_\beta$, and construct a block $B_0$ as before out of
$W\times[-1,1]$. Give this block {\em some} metric with these
properties: 
\begin{itemize}
\item Symmetry of gluing surfaces: Each component of $\boundary_\pm
B_0$ is isomestric to a {\em fixed} copy of $S_{0,3}$, which admits a
6-fold orientation-preserving symmetry group permuting the boundary
components.  

\item Flat annuli: All of the annuli of $\boundary B_0 \setminus
\boundary_\pm  B_0$ are {\em flat} -- that is, isometric to a circle
cross an interval. We assume that all the circles have length 1.
\end{itemize}

Now given any block $B(e)$, identify it with $B_0$ so that $e^-$ is
identified with  $\alpha$ and $e^+$ is identified with $\beta$. Pull
back the metric from $B_0$ to $B(e)$.

The symmetry properties imply that all the gluings can be done by
isometries (possibly after isotopy). Thus we obtain a metric on the
union of the blocks, and the boundary tori are all Euclidean. 

\subsubsection*{Geometry of the tubes}
For each vertex $v$ in $H$ we have the associated ``tube''
$C_v \times (s,t)$  in the complement of the blocks, 
which we call $U(v)$. The torus $\boundary U(v)$ has a {\em natural
marking} by  
a pair of curves -- the core curve $\gamma_v$ of $C_v\times\{s\}$, and
the meridian $\mu_v$ of $U(v)$. This marking allows us to record the geometry
of the torus via ``Teichm\"uller data'': $\boundary U$ is a Euclidean
torus in the metric inherited from the blocks, and there is a unique
number 
$$\omega \in \Hyp^2 = \{z\in\C: \Im z>0\} 
$$ 
such that $\boundary U(v)$ can be identified by an
orientation-preserving isometry with
the quotient $\C/(\Z+\omega\Z)$,
such that $\R$ and $\omega\R$ map to the classes of 
$\gamma_v$ and $\mu_v$, respectively.  We define $\omega_M(v) \equiv
\omega$, the {\em vertex coefficient} of $v$. 
Note that $|\omega_M(v)|$ is the length of the meridian $\mu_v$. 

We can then extend the metric on $\boundary U$ to make $U$ a
``hyperbolic tube'' as follows: Given $r>0$ and $\lambda\in\C$ with
$\Re\lambda>0$,  
let $\MT(\lambda,r)$ denote the quotient of an $r$-neighborhood of a
geodesic $L$ in $\Hyp^3$ by a translation $\gamma$ whose axis is $L$ and whose
complex translation distance is $\lambda$. The boundary
$\boundary\MT(\lambda,r)$ is a Euclidean torus, on which there is a
natural marking by a representative of $\gamma$ and by a meridian.  
Hence we obtain a 
Teichm\"uller coefficient $\omega(\lambda,r)$ as above. It is a straightforward
exercise to show that, given $\omega_M(v)$ there is a unique $(\lambda,r)$
such that after identifying the markings we have $\omega(\lambda,r) =
\omega_M(v)$.   We then put a metric on $U$ by identifying it with
$\MT(\lambda,r)$. 

It is not hard to check that as $|\omega_M(v)|\to\infty$, the radius
$r$ of the tube goes to $\infty$, and the length $|\lambda|$ goes to 0.

Let $M_g$ denote the union of blocks and tubes, with the metric we
have described, and the identification with $S\times\R$ we have given.

%

\section{From ending laminations to model manifold}

Given  a doubly degenerate Kleinian surface group
$\rho:\pi_1(S)\to \PSL 2(\C)$, Theorem \ref{EL exist}
gives us a pair of ending laminations $\nu_+,\nu_-$. How do these determine a
geodesic and a hierarchy from which we can build a model? 
Roughly speaking, the laminations are ``endpoints at $\infty$'' for
the hierarchy.

\subsection*{Background}

\subsubsection*{Hyperbolicity}
With Masur in \cite{masur-minsky:complex1}, we proved that
\begin{theorem}{hyperbolicity}
$\CC(S)$ is a $\delta$-hyperbolic metric space. 
\end{theorem}

We recall the definition, due to Cannon and
Gromov \cite{gromov:hypgroups,cannon:negative}: A geodesic metric
space $S$ is $\delta$-hyperbolic if all triangles are
``$\delta$-thin''. That is, given a geodesic triangle 
$[xy]\union[yz]\union[xy]$, each side is contained in a
$\delta$-neighborhood of the union of the other two. 

This simple synthetic property has many important consequences, and
gives $X$ large-scale properties analogous to those of the classical
hyperbolic space $\Hyp^n$, and any infinite metric tree. In
particular, $X$ has a {\em boundary at infinity}, $\boundary X$,
defined roughly as follows: we fix a basepoint $x_0$ 
and endow $X$ with a ``contracted'' metric $d_0$
in which $x,y\in X$ are close if
\begin{itemize}
\item they are close in the original metric of $X$, or
\item they are ``visually close'' as seen from $x_0$ -- that is, 
geodesic segments $[x_0x]$ and $[x_0y]$ have large initial segments
$[x_0x']$ and $[x_0y']$ which are in $\delta$-neighborhoods of each
other (figure \ref{visclose}). 
\end{itemize}

\realfig{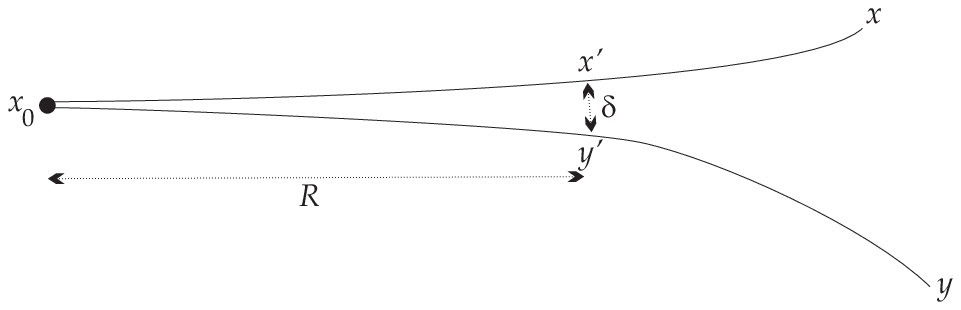}{}

The completion of $X$ in this contracted metric yields new points,
which comprise $\boundary X$. 
The construction does not in fact depend  on the choice of $x_0$. See
\cite{short:notes} for more details. The boundary of $\CC(S)$ turns
out to be a certain {\em lamination space}:

\subsubsection*{Laminations}
Thurston introduced the space of {\em measured geodesic laminations}
on a surface $S$, $\ML(S)$. Fixing a complete finite-area hyperbolic
metric on $int(S)$, a geodesic lamination is a closed subset foliated
by geodesics. A {\em transverse measure} on a geodesic lamination is a
family of Borel measures on arcs transverse to the lamination,
invariant by holonomy; that is, by sliding along the leaves. (See
Figure \ref{lamination}). 

\realfig{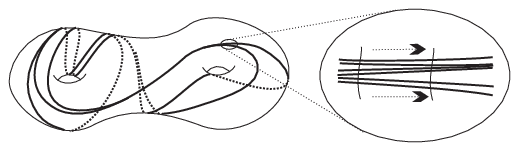}{}

This space has a natural topology, which makes $\ML(S)$ homeomorphic
to $\R^{6g-6+2n}$ when $S=S_{g,n}$. The choice of the hyperbolic
metric is not important; all choices yield naturally homeomorphic
spaces.  See Bonahon
\cite{bonahon:laminations} for more. Taking the quotient of $\ML(S)$
(minus the empty lamination) by scaling of the measures yields the
sphere $\PML(S)$ which was mentioned in Lecture 1. 
We will actually need to consider a stronger quotient,
the space of ``unmeasured laminations''
$$
\UML(S) = \ML(S)/\text{measures}
$$
Thus this is the space of all geodesic laminations which are the
supports of measures, with a quotient topology obtained by forgetting
the measure. This is different from from $\PML(S)$ because of the
existence of non-uniquely ergodic laminations, and its topology is
different from the topology on plain geodesic
laminations obtained by Hausdorff convergence of compact subsets of
$S$.  Note that the simple closed curves, i.e. vertices of $\CC(S)$,
form a dense subset of $\UML(S)$.

$\UML(S)$ is not a Hausdorff space. However, consider the
subset 
$$
\EL(S) \subset \UML(S)
$$
consisting of all ``filling'' laminations. That is, $\lambda\in\EL(S)$
if and only if all complementary regions of $\lambda$ in $S$ are ideal
polygons or once-punctured ideal polygons. An equivalent condition is
that any lamination in $\UML(S)$ different from $\lambda$ intersects
it transversely. We then have (Klarreich \cite{klarreich:boundary})
that $\EL(S)$ is a Hausdorff space. 
Furthermore, elements of $\EL(S)$ are exactly those laminations that
occur as {\em ending laminations} for manifolds $N_\rho$ without extra
parabolics. (In Theorem \ref{EL exist}, the convergence
to the ending lamination can now be understood as convergence in
$\UML(S)$.)

Klarreich showed in \cite{klarreich:boundary} that:

\begin{theorem}{C boundary}
There is a homeomorphism 
$$
k:\boundary\CC(S) \to \EL(S)
$$
such that a sequence $\beta_i\in\CC_0(S)$ converges to
$\beta\in\boundary\CC(S)$ if and only if it converges to $k(\beta)$
in the topology of $\UML(S)$.
\end{theorem}
Thus, ending laminations are points at infinity for $\CC(S)$, and from
now on we identify $\boundary \CC(S)$ with $\EL(S)$.

\subsubsection*{From lamination to hierarchy}
Now given a doubly degenerate $\rho:\pi_1(S)\to\PSL 2(\C)$, with
ending laminations $\nu_\pm$, we would like to produce a bi-infinite
geodesic $g$ in $\CC_1(S)$ whose endpoints on $\boundary\CC(S)$ are
$\nu_\pm$. 

If $\CC(S)$ were {\em locally finite}, this would be easy: Take a
sequence $\{x_i\}_{i=-\infty}^{\infty}$ in $\CC_0(S)$ such that
$$
\lim_{i\to \pm\infty} x_i  = \nu_\pm
$$
and note that, by hyperbolicity of $\CC(S)$ and the definition of
$\boundary\CC(S)$, the geodesic segments $[x_{-n},x_n]$
and $[x_{-m},x_m]$ are $2\delta$-fellow travelers on larger and larger
segments as $n,m\to\infty$.  
Thus we would expect, after extracting a subsequence, to obtain a
limiting geodesic with the endpoints $\nu_\pm$ at infinity.


For $\CC(S)$, which is {\em not} locally finite, the convergence step
is not automatic. The machinery in \cite{masur-minsky:complex2} gives
a way of getting around this, and extracting a convergent subsequence
after all. We will leave out this argument, and assume
from now on that we have a geodesic with endpoints $\nu_\pm$. 

Now, for $S=S_{0,5}$, we are ready to repeat the hierarchy (``wheel'')
construction of the previous lecture, and build from this our model
manifold. 


\medskip

Our discussion so far yields us the following: Given a
doubly-degenerate Kleinian
surface group $\rho:\pi_1(S) \to \PSL 2(\C)$, we obtain via
Bonahon-Thurston
its ending laminations $\nu_\pm\in\EL(S)$. Using Klarreich's theorem and
the work in \cite{masur-minsky:complex2}, we produce a geodesic $g$
and a hierarchy $H_\nu$, and a model manifold $\modl$ -- all depending
only on $\nu_\pm$ and not on $\rho$. Our next task will be to connect
the geometry of $\modl$ to the geometry of the hyperbolic 3-manifold
$N_\rho$.

\section{The quasiconvexity argument}

Our goal, in \S 6, is to produce a map $f:\modl \to N$ which is
uniformly Lipschitz on each of the blocks. In particular, 
if $v$ is a
vertex of $H_\nu$, it appears in some block with a {\em fixed length}
(independent of $v$, since all blocks are isometric), and so its image
has to have bounded length: 
$$
\ell_\rho(v) \le L
$$
for some uniform $L$. (Here $\ell_\rho(v)$ denotes the length in $N$
of the geodesic representative of $v$ via $\rho$).

To obtain a bound like this, we must exhibit
some connection now between the {\em geometry of $N_\rho$ } and the
{\em combinatorics/geometry of $\nu_\pm$ in $\CC(S)$}.

Recall that $\nu_\pm$ are by definition limits in $\UML(S)$ of
bounded-length curves: that is, there exists a sequence
$\{\alpha_i\}_{i=-\infty}^\infty$ in $\CC_0(S)$ with
$\ell_\rho(\alpha_i) \le L_0$ and $\lim_{i\to\pm\infty}\alpha_i =
\nu_\pm$.
The geodesic $g$ also accumulates onto $\nu_\pm$ at infinity. 
However there seems to be no a priori reason for the $\alpha_i$ to be
anywhere near $g$. Define
$$
\CC(\rho,L) = \{\alpha\in\CC_0(S): \ell_\rho(\alpha)\le L.\}
$$
To understand the relation of $g$ to $\CC(\rho,L_0)$, we will begin by
proving: 

\begin{theorem+}{Quasiconvexity}
For all $L\ge L_0$ there exists $K$, so that
$\CC(\rho,L)$ is $K$-quasiconvex.
\end{theorem+}
(Recall that $A\subset X$ is $K$-quasiconvex if for any geodesic
segment $\gamma$ with $\boundary \gamma \subset A$, $\gamma\subset
\nbhd_K(X)$.)

By hyperbolicity of $\CC(S)$ and the definition of the boundary it is
not hard to see that, since $\alpha_i$ converge to the endpoints of
$g$ as $i\to \pm\infty$,
each finite segment $G$ of $g$
is, for large enough $i$, in a
$2\delta$-neighborhood of 
$[\alpha_{-i},\alpha_i]$. 


Now $\alpha_i\in\CC(\rho,L_0)$, so using the quasiconvexity theorem,
this means that all of $g$ is in a $K'$-neighborhood of
$\CC(\rho,L_0)$.

We remark that, since $\CC$ is locally infinite, a distance bound like
this is only a weak sort of control. The Projection Bound Theorem in
Lecture 5 will be a considerably stronger generalization.

\subsection*{The bounded-curve projection}

Our main tool will be a ``coarsely defined map'' from $\CC(S)$ to
$\CC(\rho,L)$: 
$$
\Pi_{\rho,L}  : \CC(S) \to \PP(\CC(\rho,L))
$$
where $\PP(X)$ is the set of subsets of $X$. $\Pi_{\rho,L}$ ($\Pi$ for
short) will have the following properties: 
\begin{enumerate} 
\item Coarse Lipschitz: 
$$d(x,y) \le 1 \implies \diam(\Pi(x)\union\Pi(y)) \le A$$
\item Coarse Idempotence:
$$x\in\CC(\rho,L) \implies x\in\Pi(x)$$
\end{enumerate} 
(where $A$ is a constant independent of $\rho$)

These properties imply that $\Pi$ is, in a coarse sense, like a
Lipschitz projection to the set $\CC(\rho,L)$. Together with
hyperbolicity of $\CC(S)$, this has strong consequences: 

\begin{theorem}{projection implies qc}
If $X$ is $\delta$-hyperbolic, $Y\subseteq X$ and $\Pi:X\to \PP(Y)$
satisfies properties (1) and (2), then $Y$ is quasiconvex. 
\end{theorem}

The proof is similar to the proof of ``stability of quasigeodesics''
in Mostow's rigidity theorem. See \cite{minsky:boundgeom} for more
details. 

\subsection*{Definition of $\Pi_{\rho,L}$}

\subsubsection*{Pleated surfaces}

A pleated surface (or pleated map) is a map  $f:S\to N$, together with
a hyperbolic 
structure $\sigma_f$ on $S$, with the following properties: 
\begin{itemize}
\item $f$ takes $\sigma_f$-rectifiable paths in $S$ to paths in $N$ of
the same length.
\item There is a $\sigma_f$-geodesic lamination $\lambda$ on $S$, all of
whose leaves are mapped geodesically,
\item The complementary regions of $\lambda$ are mapped totally
geodesically.
\end{itemize}
We call $\sigma_f$ the {\em induced metric}, since it is determined
uniquely by the map and the first condition. The minimal $\lambda$
that works in the definition is called the {\em pleating locus} of
$f$. Informally one can think of the map as ``bent'' along $\lambda$.

This definition is due to Thurston and plays an important role in the
synthetic geometry of hyperbolic 3-manifolds. A standard example, which
we will be making use of, is the ``spun triangulation'': 

Begin with any set $P$ of curves cutting $S$ into pairs of pants, and
fix a hyperbolic metric $\sigma$ on $int(S)$ of finite area, so that
the ends are cusps.
On each component of  $P$ place one vertex, and then
triangulate each pants using only arcs terminating in these vertices
and in the cusps. 
Now ``spin'' this triangulation around $P$, by
applying a sequence of Dehn twists around each component. If at each
stage the triangulation is realized by geodesics in $\sigma$,
then the geometric limit of the
sequence will be a lamination with closed leaves $P$ and a finite
number of infinite leaves that spiral on $P$ and/or exit the cusps.
(See Figure \ref{spin})

\realfig{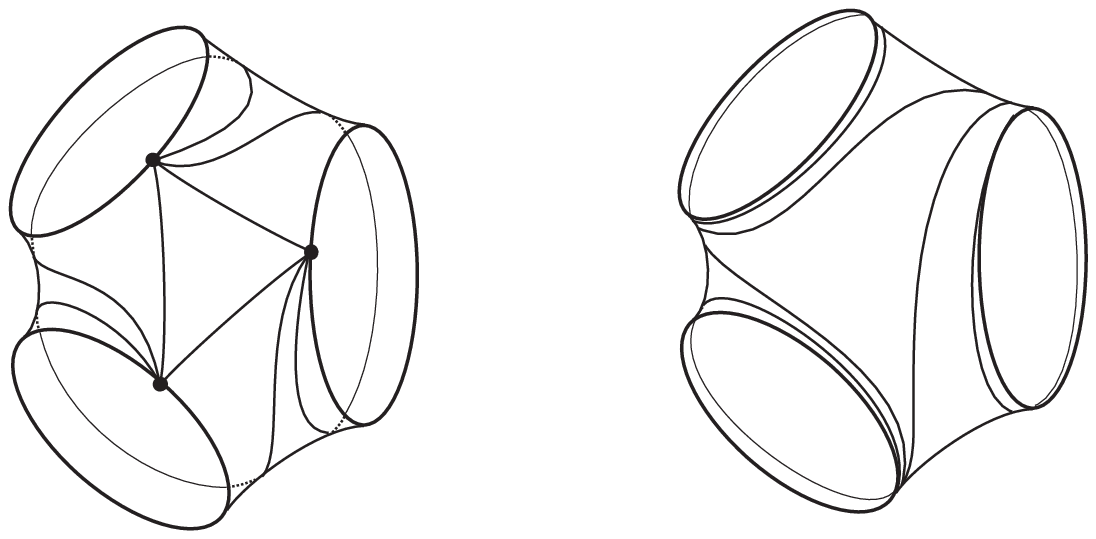}{The triangulation and ``spun'' lamination on a pair of
pants, when all boundary components are in $P$. If some are in
$\boundary S$ then the leaves go out a cusp instead of spiraling.}

In a similar way we can produce a pleated surface, first by mapping
the curves of $P$ to their geodesic representatives in $N$, and then
``spinning'' 
the images of the triangulation leaves. Finally when the leaves are in
place we fill in the spaces between them with (immersed) totally geodesic ideal
triangles, and obtain a surface together with induced metric.
(This construction is easier to visualize equivariantly in the
universal cover).

It is clear from this example that for any essential curve $\gamma$ in $S$
there is a pleated map in the homotopy class of $\rho$ that maps
$\gamma$ to its geodesic representative in $N$. We define
$$
\pleat_\rho(\gamma)
$$
to be the set of all such pleated maps. 

Now for a complete hyperbolic metric $\sigma$ on $int(S)$, define
$$
\short_L(\sigma) = \{v\in\CC_0(S): \ell_\sigma(v) \le L\}.
$$

We can now define:
\begin{equation}
\Pi_{\rho,L}(\alpha) = \bigcup_{f\in\pleat_\rho(\alpha)}
\short_L(\sigma_f).
\end{equation}

It is an observation originally of Bers that given $S$ there is a
number $L_0$ so that, for every hyperbolic metric $\sigma$ on $S$
there is a pants decomposition made up of curves of length at most
$L_0$. We call this number the ``Bers constant''. Hence for $L\ge L_0$,
$\Pi_{\rho,L}(\alpha)$ is always non-empty, and moreover
contains a pants decomposition.

Note that if $v\in\CC(\rho,L)$ then
$v\in\Pi_{\rho,L}(v)$, since if $f\in\pleat_\rho(v)$,
$\ell_{\sigma_f}(v) = \ell_\rho(v) \le L$. Hence property (2)
(Coarse Idempotence) is established.

Now our main claim, the Coarse Lipschitz property (1), will 
follow from the apparently weaker claim: 
\begin{equation}\label{bounded image}
\diam_{\CC(S)} \left(\Pi_{\rho,L}(v)\right) \le b
\end{equation}
for a priori $b$ (depending on $L$) and any simplex $v$. 

\subsubsection*{Proof of inequality (\ref{bounded image})}
First note that, for any $\sigma$,
\begin{equation}\label{bounded shorts diameter}
\diam_{\CC(S)}(\short_L(\sigma)) \le C(L).
\end{equation}
This is easy: If two curves have a length bound with respect to the
same metric $\sigma$, their intersection number is bounded in terms of
this, and a bound on the intersection number implies a bound on the
$\CC(S)$-distance by an inductive argument (see
\cite{masur-minsky:complex1}, or Hempel \cite{hempel:complex}).

Thus our main point will be to show that, for some a priori constant
$L_1$, 
\begin{equation}\label{nonempty shorts intersection}
\short_{L_1}(\sigma_f) \intersect \short_{L_1}(\sigma_g) \ne \emptyset
\end{equation}
for any $f,g\in\pleat_\rho(v)$.
This would imply, together with (\ref{bounded shorts diameter}), that
$\diam_{\CC(S)}(\Pi_{\rho,L_1}(v))\le 2C(L_1)$. In fact since
$\short_L$ is increasing with $L$ we can conclude that
$\diam_{\CC(S)}(\Pi_{\rho,L}(v))\le 2C(\max(L,L_1))$ for any $L$.

To prove inequality (\ref{nonempty shorts intersection}), let us
construct a curve 
$\gamma$ which has bounded length in both $\sigma_f$ and
$\sigma_g$. At first we note that $f(S)$ and $g(S)$ are only
guaranteed to agree on the curve $v$ itself, and this curve may be
very long. Consider the geodesic representing $v$ in the $\sigma_f$
metric on $S$. 

Suppose first that there are {\em no thin parts} on $\sigma_f$ -- that
is that $inj(\sigma_f) > \ep$ for some fixed $\ep>0$. This means that
$S$ is a closed surface, and that $\sigma_f$ has no closed geodesics
of length less than $2\ep$.
In this case we can 
approximate $v$ with a curve of bounded length that is composed of a
segment of $v$ concatenated with a very short ``jump''. That is,
for $\ep'<\ep/2$, consider the $\ep'$-neighborhood of a long segment $s$
of $v$. If this is an embedded rectangle in $S$ then its area is
at least $2\ep' l(s)$. Thus the finiteness of the area of $S$ implies
that for a certain $l(s)$, this neighborhood (which locally looks like
an embedded rectangle because of the injectivity radius lower
bound) will fail to embed globally, and where this happens we get a
``short cut'' of length $2\ep'$ joining a long (but bounded) segment
of $v$ to itself (Figure \ref{shortcut}). If we are slightly more
careful we can arrange for 
the resulting closed curve $\beta$ to be simple, and homotopically essential.

\realfig{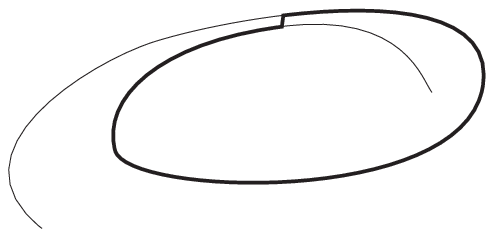}{}

Now to bound the length of $\beta$ with respect to $\sigma_g$ requires 
a bound on the $\sigma_g$ length of (the homotopy class rel endpoints
of) the short jump part of $\beta$ -- the part that runs along $v$ is
already the same length in both metrics. In other words we need to
prevent a certain kind of ``folding'' of $g$, as suggested by Figure
\ref{badfold}.

\realfig{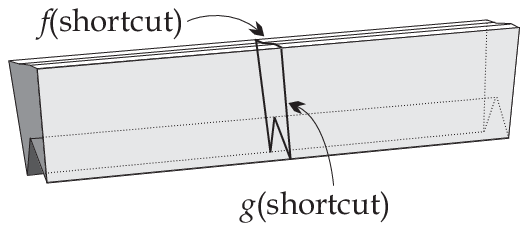}{The thin rectangle on the top is in the image of
$f$, whereas the image of $g$ is folded so that the shortcut is not
short in $\sigma_g$. This is what Uniform Injectivity rules out.}

This is prevented by a result of Thurston called the Uniform
Injectivity Theorem \cite{wpt:I}. This theorem states, in our setting,
that two leaves mapped geodesically by a pleated surface cannot line
up too closely in the image unless they are already close in the
domain. The two endpoints of the short cut part of $\beta$ are the
midpoints of long subsegments of $v$ that are close to each other in
$\sigma_f$ (we can force them to be as long as we like by taking
$\ep'$ small enough) and hence the same is true of their images by
$f$, so that the leaves line up nearly parallel in the image at
the endpoints of the short cut.
Since $f=g$ on $v$, we can then apply the Uniform Injectivity Theorem
to $g$ and conclude that the endpoints of the short cut are close
together in $\sigma_g$, and in fact (with a bit more care) the short
cut itself is homotopic rel endpoints to an arc of length at most
$\delta(\ep')$.

When $S$ has nonempty boundary,
we must take a bit more care that the closed curve
$\beta$ is in a non-peripheral homotopy class -- we may have to use
two segments on $v$ and two short cuts. 

If we allow $\sigma_f$ to have very short geodesics, we must consider
one more case. If $v$ does not enter any thin part of $S$ with core
length less than $\ep$, then the previous argument applies. If $v$
does enter the $\ep$-thin part of a $\sigma_f$-geodesic $\beta$, then
the approximation by a 
bounded-length curve may fail. However, in this case we see that both
$f$ and $g$, since they agree on $v$, have images that meet the
$\ep$-Margulis tube associated with $\beta$ in $N$, and by standard
properties of pleated surfaces this implies that $\beta$ itself has
uniformly bounded length in $\sigma_g$ as well as $\sigma_f$. 

\subsubsection*{Inequality (\ref{bounded image}) implies property (1)}
Now we are ready to establish the coarse Lipschitz property (1).
If $d(x,y) \le 1$ then $x$ and $y$ represent disjoint curves (assume
here that $S$ is not a one-holed torus or 4-holed sphere -- for those
cases there is a very similar argument). Thus the simplex $[xy]$
represents a curve system on $S$, and $\pleat_\rho([xy])$ is a
nonempty set of pleated surfaces. But it is clear that
$$
\pleat_\rho([xy]) = \pleat_\rho(x) \intersect \pleat_\rho(y)
$$
and thus this intersection is nonempty. It follows immediately that
$$
\Pi_{\rho,L}(x) \intersect \Pi_{\rho,L}(y) \ne \emptyset
$$
for any $L\ge L_0$. Hence the diameter bound (\ref{bounded image}) on 
$\Pi(x)$ and $\Pi(y)$ implies a bound on the union.

This concludes our sketch of the proof of the coarse Lipschitz
property for $\Pi_{\rho,L}$ and hence, via Theorem \ref{projection
implies qc}, of the Quasiconvexity Theorem \ref{Quasiconvexity}.

\section{Quasiconvexity and projection bounds}

The quasiconvexity of $\CC(\rho,L)$ implies that the geodesic $g$
connecting $\nu_\pm$ is a bounded distance from $\CC(\rho,L)$ (as we
discussed in lecture 4), and furthermore that
\begin{equation}\label{projection bound}
d_{\CC(S)}(v,\Pi_{\rho,L}(v)) \le B
\end{equation}
for a priori $B$ and all $v$ in $g$. However we might now wonder
{\em what good is such an estimate}, since $\CC(S)$ is locally
infinite?

Using a generalization of the Quasiconvexity Theorem, we will obtain 
a strengthening of the bound (\ref{projection bound}), which will 
then enable us in \S 6 
to establish the A Priori Bounds Theorem 
and the Lipschitz Model Theorem.

%

\subsection*{Relative bounds for subsurfaces}
In order to state our generalization of the projection bound
(\ref{projection bound}), we  must consider subsurfaces of $S$ and
their associated complexes. 

The {\em arc complex} $\AAA(W)$ of a (non-annular) surface with
boundary $W$ is  defined as follows: Vertices of $\AAA(S)$ are
homotopy classes of either essential simple closed curves (as for
$\CC(S)$) or properly embedded arcs. In the latter case the homotopy
is taken to keep the endpoints in $\boundary W$. Simplices correspond
to disjoint collections of arcs or curves.
Hence $\CC_0(W)\subset \AAA_0(W)$, and (except in the sporadic cases
$S_{0,4}$ and $S_{1,1}$, which require a separate discussion)
$\CC(W)\subset\AAA(W)$.

We also note that $\CC_0(W)$ is {\em cobounded} in $\AAA(W)$, that is,
every point of $\AAA(W)$ is a bounded distance from $\CC_0(W)$, and
that distance in $\AAA(W)$ is estimated by distance in $\CC(W)$.
Thus the two complexes are quasi-isometric.  This is easy to see with
a picture (Figure \ref{arccurve}).

\realfig{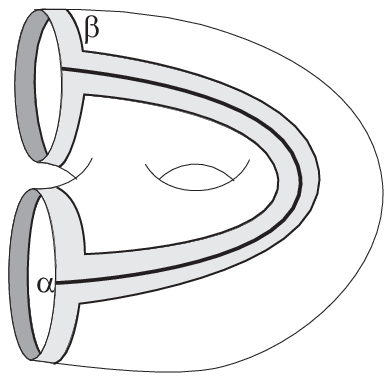}{The regular neighborhood of an arc $\alpha$ and
$\boundary W$ contains an essential curve $\beta$ in its
boundary. This gives a 
quasi-isometry from $\AAA_0(W)$ to $\CC_0(W)$.}

If $W\subset S$ is an essential subsurface, we obtain a map
$$
\pi_W : \AAA(S) \to \AAA(W)\union \{\emptyset\}
$$
defined by taking a curve system $v$ to the (barycenter of the)
simplex formed by the essential intersections $[v\intersect W]$, or to
$\emptyset$ if there are no essential intersections.

\realfig{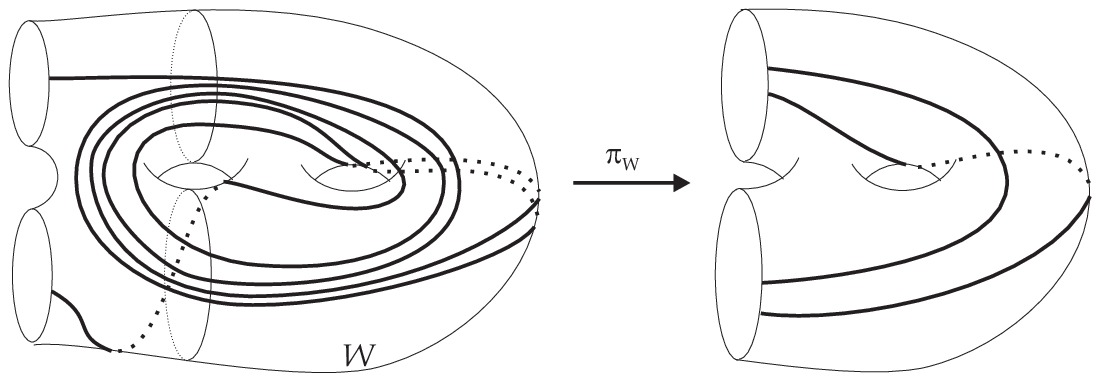}{}

For annuli in $S$ we need a different definition. Let $W\subset S$ be an
essential, nonperipheral annulus, let $\hhat W$ be the associated
annular cover of $S$, and $\bbar W$ its natural compactification.
(See Figure \ref{anncover}).

\realfig{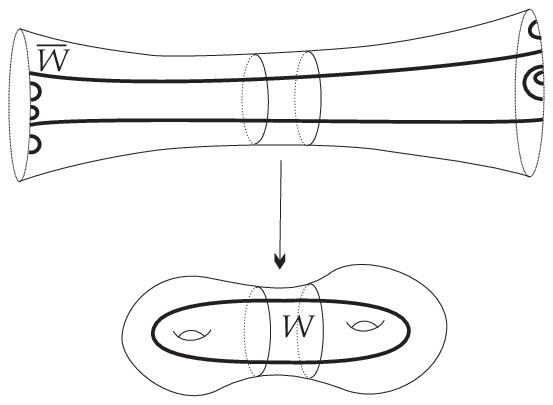}{}

Let $\AAA(\bbar W)$ be as above, except that vertices are now properly
embedded arcs up to homotopy with {\em fixed endpoints}.

Now, any $\alpha\in\AAA_0(S)$ lifts to an arc system in $\hhat W$,
which compactify to arcs in $\bbar W$. This system contains essential
arcs in $\AAA(\bbar W)$ (those with endpoints on both boundaries)
exactly if $\alpha$ intersects $W$ essentially. The set of essential
lifted arcs gives us $\pi_W(\alpha)$ in the annulus case. 

Let 
$$
d_W(\alpha,\beta) = dist_{\AAA(W)}(\pi_W(\alpha),\pi_W(\beta))
$$
(replacing $\AAA(W)$ with $\AAA(\bbar W)$ in the annulus case).
This makes sense provided both projections are nonempty. 
In the annulus case, this distance measures ``relative twisting'' of
$\alpha$ and $\beta$ around $W$.  We similarly define
$\diam_W(X) = \diam_{\AAA(W)}(\pi_W(X))$.


We can now state the
generalization of (\ref{projection bound}): 

\state{Projection Bound Theorem.}{If $v$ is in the hierarchy $H_\nu$
and $W$ is an essential subsurface other than a three-holed sphere,
then
$$
\diam_W(v\union \Pi_{\rho,L}(v)) \le B
$$
provided $\pi_W(v)$ and $\pi_W(\Pi_{\rho,L}(v))$ are nonempty,
where $L\ge L_0$ and $B$ depends on $ S$ and $L$.
}

Note that $\pi_W(\Pi_{\rho,L}(v))$ is {\em always} nonempty if $L\ge
L_0$ and $W$ is not an annulus, since then 
$\Pi_{\rho,L}(v)$ contains a pants decomposition. 

The proof of this theorem is a fancier version of the argument we used
for the Quasiconvexity Theorem in the previous lecture. 
An important ingredient is an adaptation of the ``short cut''
construction that yielded the curve $\beta$ of bounded length on
two pleated surfaces $f,g \in \pleat_\rho(v)$ (Figure
\ref{shortcut}). In the  
context of this theorem we need to make sure that $\beta$ has
essential intersection with the given subsurface $W$, and so the choice
of $\beta$ has to be carefully guided using Thurston's ``train
tracks''. In addition to this, there is an inductive structure to the
argument, using the hierarchy $H_\nu$.

\subsection*{Penetration in Margulis tubes}

We can apply the Projection Bound Theorem to control the way in which
a pleated surface enters a Margulis tube in $N$. This will then 
play an important role in the A Priori Bound Theorem in \S 6.

Let $\MT_\ep(\alpha)$ denote the {\em $\ep$-Margulis tube } in $N_\rho$
of $\rho(\alpha)$, for an element $\alpha$ of $\pi_1(S)$ (or a vertex 
$\alpha$ of $\CC(S)$). This 
is the locus where the translation length of $\rho(\alpha)$ or some power of
$\alpha$ is bounded by $\ep$. If $\ep$
is less than the Margulis constant $\ep_0$, and
$0<\ell_\rho(\alpha)<\ep$, then $\MT_\ep(\alpha)$ is a solid torus,
isometric to the hyperbolic tube $\MT(r,\lambda)$ (see \S 2)
where $\lambda$ is the
complex translation length of $\rho(\alpha)$ and 
the radius $r$ goes to $\infty$ as $\frac{\ell_\rho(\alpha)}{\ep} \to
0$. Our next goal is to detect the presence of Margulis tubes in $N$,
from the structure of the hierarchy.

\state{Tube Penetration Theorem}{ {\rm (stated for $S=S_{0,5}$)}
There exists $\ep>0$ depending on $S$, such that the following holds. 

Let $s$ be a ``spoke'' of the hierarchy $H_\nu$. If
$f\in\pleat_\rho(s)$, then
$$
f(S) \intersect \MT_\ep(\alpha) \ne \emptyset
$$
only if $\alpha$ is one of the vertices of $s$.
}
That is, the only way for $f$ to penetrate deeply into a tube is the
``obvious'' way -- by pleating along the core curve of the tube.

\subsection*{Proof of the tube penetration theorem}

We begin with this standard property of pleated surfaces (observed by
Thurston in \cite{wpt:I}): There exists $\ep_1>0$ such that, if a
pleated surface $f$ in the homotopy class of $\rho$ 
meets $\MT_{\ep_1}(\alpha)$ then $f^{-1}(\MT_{\ep_1}(\alpha))$ must be
contained in an $\ep_0$-Margulis tube in the metric $\sigma_f$ --
that is, only the thin part of $S$ is mapped into the thin part of
$N$. In particular it follows that $\ell_{\sigma_f}(\alpha) \le
\ep_0$. 

Now assume $\ep << \ep_1$.  Suppose $v$ is a vertex in $g$ crossing
$\alpha$ essentially, and $f\in\pleat_\rho(v) $ meets
$\MT_\ep(\alpha)$. Let $v^*$ denote the geodesic representative of $v$
in $N$, which is in the image of $f$.  Since $v$ crosses $\alpha$
essentially, must cross the
$\ep_0$-collar of $\alpha$ and hence by the previous paragraph
$v^*$ must meet an $\ep_0$-neighborhood of 
$\MT_\ep(\alpha)$.

In either the forward or backward direction in $g$ (suppose forward),
all vertices of $g$ after $v$ cross $\alpha$ essentially. Number them
$v=v_j, v_{j+1},\ldots$.
Eventually, $v_i^*$ for some $i>j$ is outside
$\MT_{\ep_1}(v)$, since $v_i\to \nu_+$.
Let us try to see when this happens.

\subsubsection*{Lower bound:}
Let $f\in\pleat_\rho([v_i v_{i+1}])$. If $f(S)$ meets
$\MT_{\ep_1}(\alpha)$, then both $v^*_i$ and $v^*_{i+1}$ cross through
the $\ep_0$-thin part of $S$ in the metric $\sigma_f$. Any point in
$v_i^*\intersect \MT_{\ep_1}(\alpha)$ can be connected, via an arc
in $f(S)$ of length bounded by $\ep_0$, to $v^*_{i+1}\intersect
\MT_{\ep_1}(\alpha)$. 

\realfig{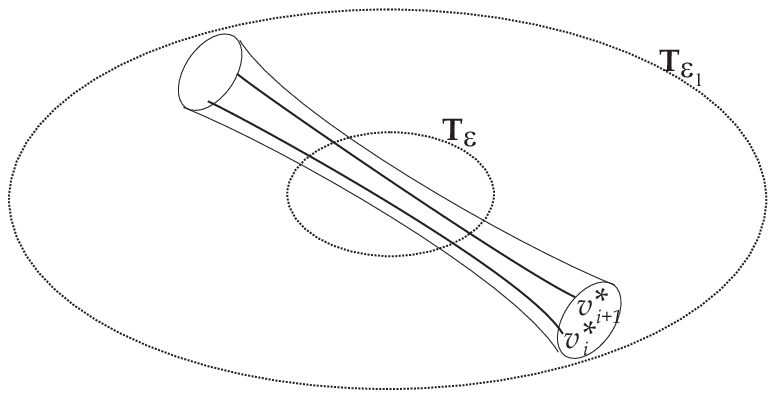}{}

Suppose then that $j+Q$ is the first value of $i$ for which $v^*_i$
fails to meet $\MT_{\ep_1}(\alpha)$. Then applying the previous
paragraph $Q$ times we have
$$
\dist(\MT_\ep(\alpha),\boundary\MT_{\ep_1}(\alpha)) \le Q\ep_0.
$$
Since, by the collar lemmas of Brooks-Matelski \cite{brooks-matelski}
and Meyerhoff \cite{meyerhoff:volumes}, this distance is an increasing
function of $\ep_1/\ep$, of the form
$$
\dist(\MT_\ep(\alpha),\boundary\MT_{\ep_1}(\alpha)) \ge
\half\log\frac{\ep_1}{\ep} - C.
$$
this gives us a lower bound  of the form
\begin{equation}\label{Q lower bound}
Q \ge a\log \frac{\ep_1}{\ep} - b.
\end{equation}

\subsubsection*{Upper bound:}
If $f\in\pleat(v_i)$ meets $\MT_{\ep_1}(\alpha)$ then, since
$\ell_{\sigma_f}(\alpha) \le \ep_0 < L_0$, we have
$$
	\alpha\in\Pi_{\rho,L_0}(v_i).
$$
The Projection Bound Theorem then implies
$$
	d_{\CC(S)}(v_i,\alpha) \le B
$$
So all such $v_i$'s lie in a ball of radius $B$. Since $g$ is a
geodesic, this means that 
\begin{equation}\label{Q upper bound}
Q \le 2B.
\end{equation}

Putting the upper and lower bounds (\ref{Q lower bound},\ref{Q upper
bound}) together, we obtain an inequality
$$
\ep > \ep_2
$$ 
where $\ep_2$ depends on the previous constants. 
Thus let us assume now that $\ep\le \ep_2$. 
Thus if $f\in\pleat_\rho(v)$ meets $\MT_\ep(\alpha)$
for $v\in g$, then $v$ and $\alpha$ must {\em not} intersect essentially.
If $v=\alpha$, we are done.

If $v\ne\alpha$ then $\alpha\in\CC(W_v)$, and we consider
the spokes $\{s_j = [u_jv]\}_{j=0}^m$ around $v$, and try to mimick the same
argument. Suppose that $u^*_j$ meets $\MT_\ep(\alpha)$, but that
$u_j \ne \alpha$.
Now $\alpha$
can be equal to at most one of the $u_i$ so let us assume it occurs
for $i<j$ if at all. Then the last vertex $u_m$ crosses $\alpha$ and,
since it is just the successor of $v$ in $g$, the previous argument applies
to it and $u_m^*$ must be outside $\MT_{\ep_2}(\alpha)$.

Now choose $Q$ to be the first positive number such that $u^*_{j+Q}$
is outside $\MT_{\ep_2}(\alpha)$. An upper bound of the form of
(\ref{Q upper bound}) follows using the same argument as before,
but applying the relative version of the Projection Bound Theorem,
$$
\diam_{W_v}(u_i,\Pi_{\rho,L_0}(u_i)) \le B.
$$

To obtain a lower bound of the form (\ref{Q lower bound}), but with
$\ep_2$ replacing $\ep_1$,  we need a
construction to replace $\pleat_\rho([v_iv_{i+1}])$, since $u_i$ and
$u_{i+1}$ do not represent disjoint curves. 

Let $\lambda_i$ denote the ``spun'' lamination, as in Lecture 4, whose
closed curves are $v$ and $u_i$. Let $\lambda_{i+1/2}$ be a ``halfway
lamination'', defined as follows. $\lambda_{i+1/2}$ contains the curve
$v$, and agrees with $\lambda_i$ and $\lambda_{i+1}$ on the complement
of $W_v$. On $W_v$ itself, $\lambda_{i+1/2}$ is as in Figure
\ref{halflam}.

\realfig{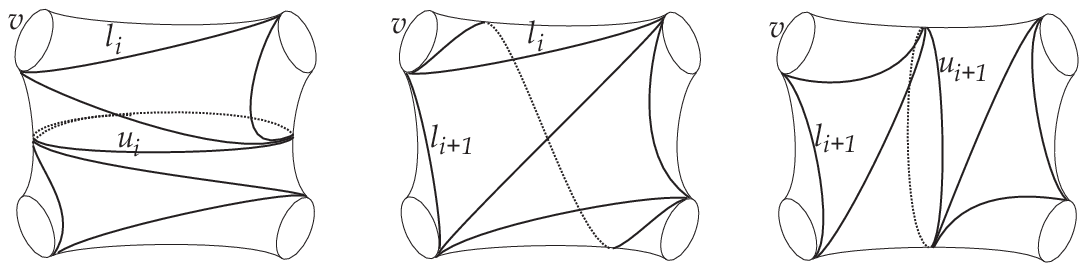}{The related laminations $\lambda_i$,
$\lambda_{i+1/2}$ and $\lambda_{i+1}$, restricted to $W_v$.}

Let $f_x$ be the pleated surface mapping $\lambda_x$ geodesically, for
$x=j+k/2$ ($k=0,\ldots , 2Q$). Note that $\lambda_{i+1/2}$ has two leaves
$l_i$ and $l_{i+1}$, which 
must cross $\alpha$ essentially since $u_i$ and $u_{i+1}$ do. 
$l_i$ is mapped to the same geodesic in $N$ by $f_{i-1/2}, f_i$ and
$f_{i+1/2}$ -- call this geodesic
$l_i^*$. We can now repeat
the lower bound argument for $Q$, finding a sequence of jumps from $u_j^*$ to
$u^*_{j+Q}$ passing through all the $l_i^*$. 
We obtain, as before, an inequality of the form 
$$
\ep > \ep_3
$$
where $\ep_3$ depends on the previous constants (and on $\ep_2$).
Thus if we choose $\ep\le \ep_3$ we must have $\alpha=u_j$ after all,
and the Tube Penetration Theorem follows. 

\medskip

For general $S$ there is an inductive argument using
the structure of the hierarchy, where the halfway surfaces need be
used only at the last stage.
\section{A-priori length bounds and model map}

In this last lecture we will sketch the proof of this basic bound: 
\state{A Priori Bound Theorem.}{%
If 
$\rho:\pi_1(S)\to \PSL 2(\C)$ is a doubly degenerate Kleinian surface
and $v$ is a vertex in the associated hierarchy $H_{\nu(\rho)}$, then
\begin{equation}\label{a priori bound}
\ell_\rho(v) \le B
\end{equation}
where $B$ depends only on the surface $S$. 
}

We will then state the Lipschitz Model Theorem and indicate how the a
priori bound is used in its proof.

\subsubsection*{Markings and elementary moves}
A marking of $S$ is a pants decomposition $\{u_i\}$ together with, for
each $i$, a transversal curve $t_i$ that is disjoint from $u_j, j\ne
i$, and intersects $u_i$ in the minimal possible way. For $S_{0,5}$, a
marking consists of 4 curves, and $t_i$ intersects $u_i$ twice. 

\realfig{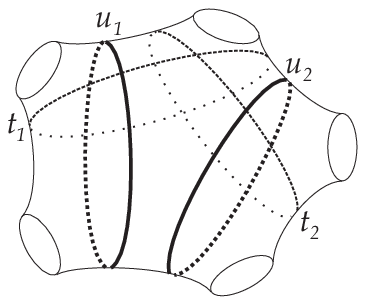}{}

An {\em elementary move} $\mu \to \mu'$ taking one marking to another
is one of the following operations: 

${\rm Twist}_i$ performs one half-twist on $t_i$ around $u_i$ (Figure
\ref{twistmove}). 

\realfig{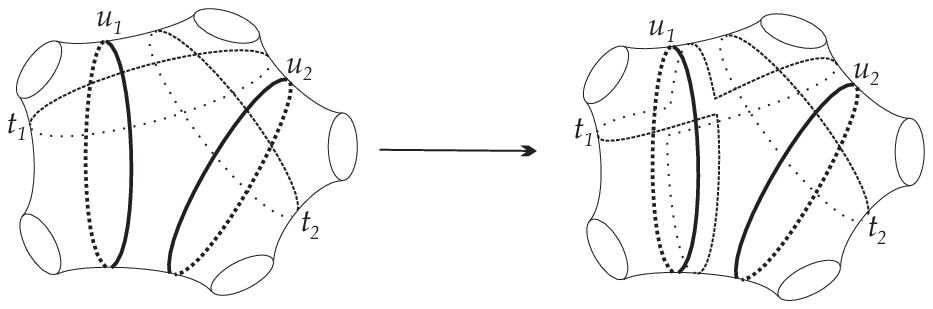}{The move ${\rm Twist}_1$}

${\rm Flip}_i$ reverses the roles of $u_i$ and $t_i$ (Figure
\ref{flipmove}). Note that in this case there has to be an adjustment
of the other $t_j, j\ne i$, so that they do not intersect the new
$u'_i$ which is the old $t_i$. There is a finite number of
``simplest'' ways to do this, and we just pick one. 

\realfig{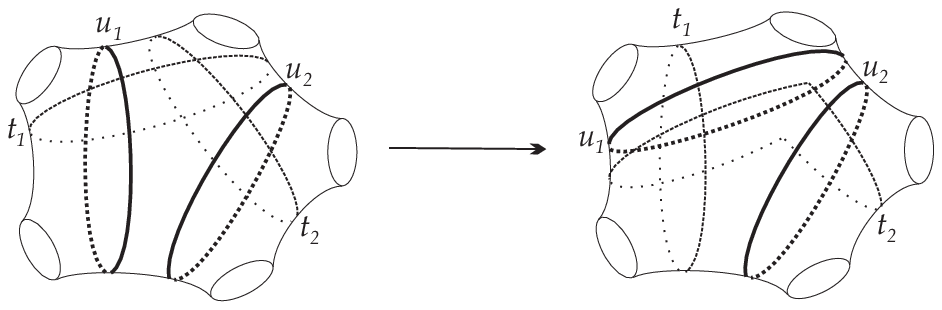}{The move ${\rm Flip}_1$}

The graph whose vertices are markings and whose edges are elementary
moves is connected, locally finite, and its quotient by the mapping
class group 
of $S$ is finite. The total length of a marking in a hyperbolic metric
$\sigma$ is just the sum of the lengths of the curves $u_i$ and
$t_i$. Note that a bound on the total length of $\mu$ in $\sigma$
constrains $\sigma$ to a bounded subset of Teichm\"uller space. 
We will need the following observation: 
If $\mu_0$ is a marking of total length $L$ in a hyperbolic metric
$\sigma$ on $S$, and 
$$
\mu_0 \to \mu_1 \to \cdots \to \mu_n
$$
is a sequence of elementary moves, then $\mu_n$ has total length at
most $K$, where $K$ depends only on $L$ and $n$.

We will control elementary moves using this theorem: 

\begin{theorem}{bound moves}{\rm (Masur-Minsky
\cite{masur-minsky:complex2})}
If $\mu$ and $\mu'$ are two markings and 
$$
\sup_{W\subseteq S} d_W(\mu,\mu') \le M.
$$
Then  there exists a sequence of elementary moves from $\mu$ to $\mu'$
with at most 
$$
CM^\xi
$$
steps, where $C$ and $\xi$ depend only on $S$.
\end{theorem}

Here, $d_W(\mu,\mu')$ is defined as in \S5, with the projection
$\pi_W(\mu)$ simply being the union of $\pi_W(a)$ over components $a$
of $\mu$. The supremum is over all essential subsurfaces in $S$,
including $S$ itself.
The proof of this theorem uses the hierarchy machinery discussed
in \S2.

\subsection*{Proving the A Priori Bounds}
Again, we are working in the case that $S=S_{0,5}$. Let $s=[u_1u_2]$
be a spoke of the hierarchy, and $f\in\pleat_\rho(s)$.
Let $\ep>0$ be the constant given by the Tube Penetration Theorem.

\subsubsection*{Thick case}
Suppose $\sigma_f$ has no geodesics of length $\ep$ or less. 
Then $\sigma_f$ admits a marking $\mu$ of total length at most $L$
(depending on $\ep$). 


The Projection Bound Theorem implies that 
$$
d_W([u_1 u_2], \mu) \le B
$$
For any $W$ meeting $u_1$ or $u_2$ essentially, since 
$\mu\subset\Pi_{\rho,L}([u_1u_2])$.
This includes all $W$ except for the annuli $A_i$ with cores $u_i$
($i=1,2$). We therefore choose transversals
$t_1$, $t_2$ for $u_1,u_2$, such that 
$d_{A_i}(t_i,\mu')$ is at most 2. Thus $\{u_1,t_1,u_2,t_2\}$ give us
a marking $\mu'$ such that
$$
d_W(\mu',\mu) \le B
$$
for {\em all } $W\subseteq S$. 
Applying Theorem \ref{bound moves}, we bound the elementary-move
distance from $\mu$ to $\mu'$, and hence obtain a bound on
$\ell_{\sigma_f}(\mu')$. 

This in turn bounds $\ell_\rho(u_i)$, which gives the a priori bound
(\ref{a priori bound}) in this case. 

\subsubsection*{Thin case}
Suppose that $\sigma_f$ does have some curve $\alpha$ of length less
than $\ep$.  Then $f(S)$ meets $\MT_\ep(\alpha)$. The Tube
Penetration Theorem now implies that $\alpha  = u_1$ or $\alpha=u_2$. 
Suppose the former, without loss of generality.
Thus, we repeat the argument of case (1) on the subsurface
$W_\alpha$, finding a minimal-length marking (now consisting of just
one curve and its transversal) and using the Projection Bound Theorem
to bound the length of $u_2$.  Again we have the a priori bound
(\ref{a priori bound}).

\subsection*{Constructing the Lipschitz map}
We are now ready to state, and summarize the proof of, our main theorem:

It will be convenient to define $\UU[k]$ to be the union of tubes 
$\{U(v):|\omega_M(v)| \ge k\}$, and to let $\modl[k] = \modl \setminus
\UU[k]$. Thus $\modl[0]$ is the union of blocks. 

If $v$ is a vertex of $H_\nu$, let $\MT_{\ep_0}(v)$ be the Margulis
tube (if any) of the homotopy class $\rho(v)$ in $N_\rho$. Let
$\MT[k]$ denote the union of $\MT_{\ep_0}(v)$ over all $v$ with
$|\omega_M(v)|\ge k$. 

\state{Lipschitz Model Theorem.}{%
There exist $K,k_0>0$ such that, if
$\rho:\pi_1(S)\to \PSL 2(\C)$ is a doubly degenerate Kleinian surface
group with end invariants $\nu(\rho)$,
then there is a map
$$
F: \modl \to N_\rho
$$
with the following properties:
\begin{enumerate}
\item $F$ induces $\rho$ on $\pi_1$, is proper, and has degree 1.
\item $F$ is $K$-Lipschitz on $\modl[k_0]$.
\item $F$ maps $\UU[k_0]$ to $\MT[k_0]$, and $\modl[k_0]$ to
$N_\rho\setminus \MT[k_0]$.
\end{enumerate}
}
Note that the Lipschitz property in part (2) is with respect to the
{\em path metric} on $\modl[k_0]$ (the distance function of $\modl$
restricted to $\modl[k_0]$ may be smaller).

The map $F:\modl[0] \to N$ can be constructed in each block $B_e$ individually.
We first define the map on the gluing boundaries $\boundary_\pm B$, which
are all isometric to a fixed three-holed sphere; let $Y$ 
be such a boundary component. As in Lecture 4, 
there is a pleated map $h:Y \to N$ 
in the homotopy class determined by $\rho$, which sends each 
boundary component of $Y$ to its geodesic representative or to the
corresponding cusp if it is parabolic (more
accurately $h$ is defined 
on $int(Y)$, and gives a metric whose completion has geodesic
boundaries for the non-parabolic ends of $Y$ and cusps for the
parabolic ends).

The a priori bounds give an upper bound on the boundary lengths of
this pleated surface. Thus, after excising standard collar
neighborhoods of the boundaries, we obtain a surface which can be
identified, with uniform bilipschitz distortion, with $Y$ under its
original model metric.
Composing this identification
with $h$, we obtain the map $F|_Y$.

Our next step is to define $F$ on the ``middle'' surface of a block, 
which we can write as $W\times\{0\}$ with $W$ a four-holed sphere
(for general $S$, $W$ can also be a one-holed torus). The
``halfway surfaces'' from the proof of the Tube Penetration Theorem
(defined using the two vertices $e^+$ and $e^-$)
provide us with a map $W\times\{0\}\to N$ and an induced hyperbolic
metric. Another application of the a-priori length bounds together
with Thurston's Efficiency of Pleated Surfaces \cite{wpt:II} implies
that this metric on $W$ is within uniform bilipschitz distortion of
the model metric. 

We then extend to the rest of the block by a map that takes vertical
lines to geodesics. A Lipschitz bound on this part of the map is then
an application of the ``figure-8'' argument from
\cite{minsky:torus}. (In brief, we let $X$ be a wedge 
of two circles in the gluing boundary $Y$ that generating a nonabelian
subgroup of $\pi_1(Y)$ and 
having bounded length in both $W\times\{0\}$ and $Y$. The extension
gives a map of $X\times[0,1]$ with geodesic tracks $\{x\}\times[0,1]$,
and if the track lengths are too long then the images of the two
circles in the middle are either both very short, or nearly parallel,
violating discreteness either way).

Fixing  $k_0$, the tubes $U$ with $|\omega(U)| < k_0$
fall into some finite set of isometry types, and the map can be
extended to those, again with some uniform Lipschitz bounds. 
If $|\omega(U)| > k_0$ then there is an {\em upper bound} for the
corresponding vertex, 
$\ell_\rho(v) \le \ep$ where $\ep\to 0$ as $k_0\to \infty$. This is
the main result of \cite{minsky:kgcc}, which uses similar tools but is
slightly different than what we have seen so far. Thus choosing $k_0$
appropriately, $U$ must correspond to a Margulis tube $\MT$ with very large
radius (and short core). The map on the blocks cannot penetrate more
than a bounded distance into such a tube (by the Lipschitz bounds) and
so composition with an additional retraction
on a collar of $\boundary \MT$ yields us a map
$F:\modl[k_0] \to N\setminus\MT[k_0]$ which takes each $\boundary U$ 
with $|\omega(U)|\ge k_0$ to the boundary of the corresonding $\MT$.

We fill in the map on the remaining $U$ and the cusps (with no
Lipschitz bounds) to obtain a final map $F:\modl \to N$ with all the
desired properties. The fact that $F$ is proper
essentially follows from the fact that every block contains a  bounded
length curve in a unique homotopy class, and thus their images cannot
accumulate in a compact set. That $F$ has degree 1 follows from the
fact that blocks far toward the $+$ end of the hierarchy correspond to
vertices close to $\nu_+$, and hence their images must go out the $+$
end of $N$.

\subsection*{Consequences}
It is now easy to obtain the following lower bound on lengths:

\begin{corollary}{short curves vertices}
There exists an a priori $\ep>0$ such that 
all curves of length less than $\ep$ in $N$ must occur as vertices of
the hierarchy.
\end{corollary}

\noindent {\em Proof:}
If $\ell_N(\gamma)<\ep$ then $\gamma$ has a Margulis tube
$\MT_\ep(\gamma)$ in $N$. Since $F$ has degree 1, this tube is in the
image of $F$. The Lipschitz bound on $F|_{\modl[k_0]}$ keeps it out of
$\MT_\ep(\gamma)$ if $\ep$ is sufficiently small. Hence $\MT_\ep(\gamma)$ is
in the image of some tube $U(v)$, which means that $\gamma$ corresponds
to $v$.

(Note that this corollary does not follow directly from the Tube
Penetration Lemma since that lemma assumes the existence of a pleated
surface associated to a spoke which penetrates $\MT_\ep(\gamma)$. The
global degree argument is necessary).

\medskip

Another corollary, which requires a bit more work and uses Otal's
theorem \cite{otal:knotting} on the unknottedness of short curves in
$N$, is the following, which describes the
topological structure of the set of short curves in $N$:

\begin{corollary}{topology of  tubes}
$\modl[k_0]$ is homeomorphic to $N\setminus \MT[k_0]$.
\end{corollary}
For a proof of this see \cite{brock-canary-minsky:ELCII}.

These results give us a complete description of the ``short curves''
in $N$, and in particular give a combinatorial criterion (in terms of
the hierarchy and the coefficients $\omega_M$) for when the manifold
has bounded geometry. This answers the short list of questions posed in the
introduction. 

We can also obtain somewhat more explicit lower bounds on $\rho$-lengths of
the vertices of the hierarchy, namely:

\begin{corollary}{upper length bounds}
If $v$ is a vertex of the hierarchy $H_{\nu(\rho)}$ then
$$
\ell_\rho(v) \ge \frac c {|\Momega(v)|^2}
$$
where $c$ depends only on the surface $S$.
\end{corollary}

This follows from the fact that, because of the Lipschitz property of
the map $F$,  $|\Momega(v)|$ bounds the meridian length of the
corresponding Margulis tube $\MT_{\ep_0}(v)$ in $N$. This gives an
upper bound for radius of the tube, and 
the collar lemmas of \cite{brooks-matelski} and
\cite{meyerhoff:volumes} then give a lower bound 
for the core length $\ell_\rho(v)$. An upper bound for 
$\ell_\rho(v)$ which goes to 0 as $|\omega_M(v)|\to \infty$ also
exists: it follows from the main theorems of \cite{minsky:kgcc}, and
was already used briefly in the proof of the Lipschitz Model Theorem.


\newcommand{\etalchar}[1]{$^{#1}$}
\ifx\undefined\bysame
\newcommand{\bysame}{\leavevmode\hbox to3em{\hrulefill}\,}
\fi

\end{document}